\documentclass[11pt]{amsart}
\usepackage{graphicx,rlepsf,amssymb,amsmath}
\usepackage[all]{xy}
\CompileMatrices
\theoremstyle{plain}
\newtheorem{theorem*}{Theorem}
\newtheorem{theorem}{Theorem}[section]
\newtheorem{lemma}{Lemma}[section]

\newtheorem{claim}{Claim}[section]
\newtheorem{corollary}{Corollary}[section]
\newtheorem{corollary*}{Corollary}
\theoremstyle{definition}
\newtheorem{definition}{Definition}[section]
\newtheorem{example}{Example}[section]
\theoremstyle{remark}
\newtheorem{remark}{Remark}[section]

\def \Z{\mathbb{Z}}
\def \Q{\mathbb{Q}}
\def \F{\mathbb{F}}
\def \R{\mathbb{R}}

\def \aug{{\rm aug}}
\def \dd{{\rm d}}

\def \E{{\rm E}}
\def \Eul{{\rm Eul}}
\def \Hom{{\rm Hom}}
\def \Id{{\rm Id}}
\def \ind{{\rm ind}}
\def \int{{\rm int}}
\def \incl{{\rm incl}}

\def \Ker{{\rm Ker}}
\def \Map{{\rm Map}}
\def \mod{{\rm mod}}

\def \Parall{{\rm Parall}}
\def \pr{{\rm pr}}
\def \rk{{\rm rk}}
\def \sgn{{\rm sgn}}
\def \T{{\rm T}}
\def \Trace{{\rm Trace}}
\def \Tors{{\rm Tors}}

\begin{document}
\title[]{Some finiteness properties for the Reidemeister--Turaev torsion of three-manifolds}

\date{February 8, 2009}

\author[]{Gw\'ena\"el Massuyeau}



\begin{abstract}
We prove for the Reidemeister--Turaev torsion of closed oriented three-mani\-folds
some finiteness properties in the sense of Goussarov and Habiro, that is,
with respect to some cut-and-paste operations which preserve the homology type of the manifolds. 
In general, those properties require the manifolds to come
equipped with an Euler structure and a homological parametrization.
\end{abstract}

\maketitle

\tableofcontents

\section{Introduction}

The theory of finite-type invariants of $3$-manifolds aims at understanding 
how manifolds are related one to the other through cut-and-paste operations
and, consequently, how their invariants behave with respect to such operations.

In the Goussarov--Habiro theory, manifolds are modified using
surgery operations which preserve the homology type \cite{GGP,Habiro,Garoufalidis}.
Given a closed oriented connected $3$-manifold $M$, a handlebody $H \subset M$
and a Torelli automorphism $h$ of $\partial H$ (that is, $h:\partial H \to \partial H$ is a diffeomorphism
which acts trivially in homology), one can form a new closed oriented connected $3$-manifold:
$$
M_h := (M \setminus \int\ H) \cup_h H.
$$
The move $M \leadsto M_h$ is called a \emph{Torelli surgery}. Thanks to Matveev's theorem \cite{Matveev},
one can always decide when two given manifolds are related by a finite sequence of such operations.
Let now $f$ be an invariant of closed oriented connected $3$-manifolds with values in $A$, an Abelian group.
It is said to be a \emph{finite-type invariant} of \emph{degree} at most $d$ if,
for any manifold $M$ and for any family $\Gamma$ of
$d+1$ Torelli automorphisms of the boundaries of pairwise disjoint handlebodies in $M$,
the following identity holds:
\begin{equation}
\label{eq:fti}
\sum_{\Gamma' \subset \Gamma} (-1)^{|\Gamma'|} \cdot f\left(M_{\Gamma'}\right) =0 \in A.
\end{equation}
Here, $M_{\Gamma'}$ denotes the manifold obtained from $M$ 
by the simultaneous surgery defined by those elements of  $\Gamma$ shortlisted in $\Gamma'$.
For integral homology $3$-spheres, this definition is essentially equivalent 
to the original notion of finite-type invariant introduced by Ohtsuki \cite{Ohtsuki}.

In the case of rational homology $3$-spheres $M$, there exist some very powerful invariants
\cite{LMO,BGRT,KT} which are universal among finite-type invariants $f$ with values in $A=\Q$. 
But, either those invariants are not defined for manifolds $M$  whose first Betti number is positive, 
either they become trivial when the latter gets too high. As a matter of fact, only very few instances of finite-type 
invariants are known for manifolds with arbitrary homology.

On the other hand, the Reidemeister--Turaev torsion  (or, \emph{RT torsion})
is quite a well-understood invariant
of a closed oriented connected $3$-manifold $M$. Denoted by
$$
\tau(M,\xi)\in Q \left( \Z [ H_1(M;\Z) ] \right),
$$
it takes its values in the quotient ring of the group ring $\Z[H_1(M;\Z)]$. 
It generalizes both the Reidemeister torsion of a lens space 
and the Alexander polynomial \cite{Turaev_first,Turaev_smooth}.
To be defined without indeterminacy, it needs the manifold $M$ to be endowed with an Euler structure $\xi$. 
Despite the combinatorial nature of the Reidemeister torsion, an Euler structure admits a geometric description: 
This is a non-singular (that is, nowhere zero) vector field on $M$ 
up to punctured homotopy (that is, up to homotopy on $M$ deprived of one point).

The question of how the RT torsion connects to finite-type invariants
seems not to have been addressed yet, 
except in the case of the one-variable Alexander polynomial of manifolds \cite{GH,CM,Lieberum}.
Nevertheless, finiteness results are quite expected for the RT torsion, 
since the multi-variable Alexander polynomial of links has finiteness properties 
with respect to the Goussarov--Vassiliev theory of finite-type invariants \cite{Murakami}.\\

To answer partly to that question, we need to abstract the homology of manifolds.
More precisely, we \emph{fix} a finitely generated Abelian group $G$, and
we consider  triples of the form
$$
(M,\xi,\psi)
$$
where $M$ is a closed oriented connected $3$-manifold, $\xi$ is an Euler structure on $M$ and 
$\psi: G \to H_1(M;\Z)$ is an isomorphism. We call $\psi$ a \emph{homological parametrization} for the closed manifold $M$,
which is like a numbering of the components for an oriented link in $S^3$.
Note that the RT torsion gives an invariant of triples $(M,\xi,\psi)$, by setting
$$
\tau(M,\xi,\psi):= Q(\psi^{-1}) \left( \tau(M,\xi) \right) \in Q(\Z[G])
$$
where $Q(\psi^{-1}): Q\left(H_1(M;\Z)\right)\to Q(\Z[G])$ is the ring isomorphism induced by $\psi^{-1}$. 

The main property of  a Torelli surgery is to preserve the homology of the manifold.
Indeed, the move $M \leadsto M_h$ comes with a canonical isomorphism 
$$
\Phi_h: H_1(M;\Z) \longrightarrow H_1(M_h;\Z),
$$ 
as it is given by the Mayer--Vietoris theorem. Another property is to define a canonical correspondence
$$
\Omega_h:\Eul(M) \longrightarrow \Eul(M_h)
$$ 
between Euler structures \cite{DM}. 
This is defined by cutting and pasting vector fields in an appropriate way 
(see \S \ref{subsec:Torelli}).
It follows that the Torelli surgery makes sense in the context of manifolds 
with Euler structure and parametrized homology:
$$
(M,\xi,\psi) \leadsto (M,\xi,\psi)_h:=\left(M_h,\Omega_h(\xi),\Phi_h \circ \psi \right).
$$
In particular, by writing $(M,\xi,\psi)$ in condition (\ref{eq:fti}) in place of $M$, one gets
the definition of a finite-type invariant for such triples 
(see \S \ref{subsec:GH_with_structures}).
The invariant $\tau(M,\xi,\psi)$, when reduced modulo a power of the augmentation ideal
$$
I := \Ker \left( \aug: \Z[G] \longrightarrow \Z \right),
$$ 
is finite-type in this sense.

\begin{theorem*}
\label{th:FTI}
Assume that $G$ has positive rank or that it is finite cyclic.
Let $d\geq 1$ be an integer. Then, the RT torsion reduced modulo $I^d$
$$
\tau(M,\xi,\psi) \in Q(\Z[G])/I^d
$$
of closed oriented connected $3$-manifolds $M$ with Euler structure $\xi$ 
and homological parametrization $\psi: G \to H_1(M;\Z)$,
is a finite-type invariant of degree at most $d+1$.
\end{theorem*}

\noindent
The proof of Theorem \ref{th:FTI} uses Heegaard splittings of manifolds. First,
any Heegaard splitting of a manifold $M$ induces a cell decomposition of it, 
with which the RT torsion can be combinatorially computed. 
This formula appears in the proof of several results by Turaev \cite{Turaev_spinc,Turaev_bigbook},
to evaluate how deep inside the $I$-adic filtration (or related filtrations) 
the RT torsion should live and to compute the leading term from the cohomology ring.
Second, following Hutchings and Lee \cite{HL}, one can pass from combinatorial Euler structures
to geometric ones using the gradient of a Morse function inducing the Heegaard splitting.
These will be our two tools to compare $\tau(M,\xi,\psi)$ to $\tau\left((M,\xi,\psi)_h\right)$,
after a Torelli surgery $(M,\xi,\psi) \leadsto (M,\xi,\psi)_h$ has been performed.
Theorem \ref{th:FTI} is deduced from a general lemma giving, for any group $G$,
finiteness properties for an invariant equivalent to $\tau(M,\xi,\psi)$.
See \S \ref{subsec:computation} and \S \ref{subsec:finiteness_properties}.

The same techniques are applied in \S \ref{subsec:non-variation}
to find sufficient conditions on the Torelli surgery instructions 
for the RT torsion not to be changed.

\begin{theorem*}
\label{th:twist}
Let $M$ be a closed oriented connected $3$-manifold with Euler structure $\xi$
and homological parametrization $\psi: G \to H_1(M;\Z)$.
Let $H$ be a handlebody in $M$ and let $h$ be a Torelli automorphism of $\partial H$
such that one of the following two conditions holds:
\begin{itemize}
\item[--] The handles of $H$ are null-homologous in $M$.
\item[--] The diffeomorphism $h$ acts trivially 
on the second solvable quotient $\pi/\pi''$ of the fundamental group $\pi$ of $\partial H$.
\end{itemize}
Then, we have that $\tau(M,\xi,\psi)= \tau\left((M,\xi,\psi)_h\right) \in Q\left(\Z[G]\right)$.
\end{theorem*}

There are alternative theories of finite-type invariants for $3$-manifolds,
all being equivalent to the Ohtsuki theory for integral homology $3$-spheres 
(up to a linear re-scaling of the degrees, and as far as rational-valued invariants are concerned \cite{GGP}). 
For instance, instead of considering any Torelli automorphisms in condition (\ref{eq:fti}),
one could be more restrictive and consider diffeomorphisms 
wich act trivially on the $c$-th nilpotent quotient $\pi/\pi_{c+1}$ of the fundamental group $\pi$ of the surface. 
For $c=1$, this is the Goussarov--Habiro theory but, for $c=2$, one gets a different theory \cite{Garoufalidis}.
Our methods give similar results of finiteness for any class $c\geq 2$, see \S \ref{subsec:Johnson}.
Nevertheless, they do not apply to the Cochran--Melvin theory \cite{CM}, 
in which finiteness properties have already been proved for the one-variable Alexander polynomial.

Theorem \ref{th:FTI} implies some well-known finiteness properties for invariants
that are known to be determined by the RT torsion. 
Those include the Casson--Walker--Lescop invariant when $\rk(G)>0$, 
and the one-variable Alexander polynomial when $\rk(G)=1$. This is checked in \S \ref{subsec:check}.

Finally, one may ask the question of whether finite-type invariants 
dominate the RT torsion. We do not have a general answer but, by an algebraic fact,
this certainly holds true when $G$ has no two elements of finite coprime orders.
In particular, the Milnor--Turaev torsion (or, \emph{MT torsion}),
which is an enhancement of the multi-variable Alexander polynomial,
is dominated by finite-type invariants. See \S \ref{subsec:domination} and \S \ref{subsec:MT}.\\

\noindent
\textbf{Acknowledgments.} The author would like to thank Vladimir Turaev 
for his suggestions and comments on this paper. 
He is indebted to the European Commission for support (MEIF-CT-2003-500246)
and the University of Pisa for hospitality, with special thanks to Carlo Petronio.

\section{A quick review of the RT torsion}
\label{sec:review_RT}

In this expository section, we briefly review the theory of RT torsion.
References on this topic include the papers \cite{Turaev_first,Turaev_smooth}
and the monographs \cite{Turaev_smallbook,Nicolaescu,Turaev_bigbook}, 
to which the reader is refered for details and proofs.
On the way, we  fix the notations that are used throughout the paper,
starting with the following conventions:
\begin{itemize}
\item[$\centerdot$] An Abelian group $G$, or its action on a set, is written additively,
except when it is seen as a subgroup of the group of units of $\Z[G]$. 
\item[$\centerdot$] Unless otherwise mentioned, (co)homology groups are computed with integral coefficients.
\end{itemize}

\subsection{RT torsion of a CW-complex}

First, one needs to define what is the torsion of a CW-complex.

\subsubsection{Reidemeister torsion of a chain complex}

Let $\F$ be a commutative field.

Given a finite-dimensional $\F$-vector field $V$ and two of its basis $b$ and $c$, 
$[b/c]\in \F\setminus \{0\}$ denotes the determinant of the matrix expressing 
$b$ in the basis $c$. The basis $b$ and $c$ are \emph{equivalent} when $[b/c]=1$.
Given a short exact sequence of $\F$-vector spaces $0 \to V' \to V \to V'' \to 0$
and basis $c'$ and $c''$ of $V'$ and $V''$ respectively, denote by $c' c''$
the equivalence class of basis of $V$ obtained by juxtaposing the image of $c'$ with a lift of $c''$.

Consider a chain complex of finite-dimensional $\F$-vector spaces
$$
C=\left( \xymatrix{ C_m\ar[r]^-{\partial_{m-1}} & C_{m-1} \ar[r] 
 & \cdots \ar[r]^{\partial_0} & C_0} \right)
$$
which comes equipped with a \emph{basis} $c$ and a
\emph{homological basis} $h$. This means that  
$c=(c_0,\dots,c_m)$ where $c_i$ is a basis of the $i$-chains space $C_i$,
and $h=(h_0,\dots,h_m)$ where $h_i$ is a basis of the $i$-th homology group $H_i(C)$.
Define the modulo $2$ integer
$$
N(C) :=\sum_{k=0}^m \alpha_k(C)\cdot \beta_k(C) \in \Z_2,
$$
$$
\hbox{where} \quad
\alpha_i(C)  := \sum_{k=0}^i \dim\left(C_k\right) \in \Z_2 
\quad \hbox{and} \quad
\beta_i(C)  :=  \sum_{k=0}^i \dim\left(H_k(C)\right) \in \Z_2.
$$
Choose, for each $i=0,\dots,m$, 
a basis $b_i$ of the space of $i$-boundaries $B_i(C)$.
Thanks to the short exact sequences
$0  \to B_i(C) \to Z_i(C) \to H_i(C) \to 0$
and $0  \to Z_i(C) \to C_i \to  B_{i-1}(C) \to 0$,
one gets a basis for $C_i$, namely $b_ih_ib_{i-1}:=(b_ih_i)b_{i-1}$.

\begin{definition} 
\label{def:algebraic_torsion}
The \emph{Reidemeister torsion} of the $\F$-complex $C$,
based by $c$ and homologically based by $h$, is
\begin{displaymath} 
\tau(C;c,h):= (-1)^{N(C)} \cdot
\prod_{i=0}^m\left[b_ih_ib_{i-1}/c_i\right]^{(-1)^{i+1}} \in \F\setminus \{0\}
\end{displaymath}
and does not depend on the choice of $b$. 
\end{definition}

\begin{remark}
This is Turaev's sign version of the Reidemeister torsion \cite{Turaev_smooth}.
\end{remark}

\subsubsection{Some structures on a CW-complex}

Let $X$ be a finite, connected CW-complex 
with Euler characteristic $\chi(X)=0$.

A \emph{homological orientation} $\omega$ of $X$ is an orientation
of the $\R$-vector space $H_*(X;\R)$, the opposite orientation
being denoted by $-\omega$.

An \emph{Euler chain} in $X$ is a singular $1$-chain $c$ on $X$ 
with boundary 
$$
\partial c = \sum_{\sigma, \textrm{ cell of } X}
(-1)^{\dim(\sigma)} \cdot c_\sigma
$$
where $c_\sigma$ denotes the center of the cell $\sigma$. 
When considered up to homology, Euler chains are called
\emph{Euler structures} and form a set $\Eul(X)$ which,
with the obvious action, is a $H_1(X)$-affine space. 

Fix a base point $\star \in X$ to determine the maximal Abelian covering 
$p: \widehat{X} \to X$. Its group of transformations is identified with $H_1(X)$. 
The cell decomposition of $X$ lifts to a cell decomposition of $\widehat{X}$. 
A family $\widehat{e}$ of cells of $\widehat{X}$ is \emph{fundamental}
when each cell $\sigma$ of $X$ has a unique lift in it, which is then denoted by $\widehat{e}(\sigma)$.
Two  fundamental families of cells $\widehat{e}$ and $\widehat{f}$ are equivalent when the difference
$$
\widehat{f}-\widehat{e}\ :=
\sum_{\sigma,\textrm{ cell of } X}
(-1)^{\dim(\sigma)} \cdot \left(\widehat{f}(\sigma)-\widehat{e}(\sigma) \right) \ \in H_1(X)
$$
vanishes. (Here $\widehat{f}(\sigma)- \widehat{e}(\sigma) \in H_1(X)$ denotes the transformation needed
to move $\widehat{e}(\sigma)$ to $\widehat{f}(\sigma)$.) When considered up to equivalence, 
fundamental families of cells form a set $\E(X)$ which is a $H_1(X)$-affine space.

Euler structures are used as ``instructions to lift cells''. Specifically, 
given a fundamental family of cells $\widehat{e}$, connect by an oriented path 
the center of each cell $\widehat{e}(\sigma)$ to a single point in $\widehat{X}$.
This path goes from $\widehat{e}(\sigma)$ to the single point 
if $\dim(\sigma)$ is odd, and vice-versa if $\dim(\sigma)$ is even. 
The image by $p$ of this $1$-chain is an Euler chain (shaped like a spider). 
Thus, one gets a $H_1(X)$-equivariant bijection $\E(X) \to \Eul(X)$.

\subsubsection{Definition}

Let $\varphi : \Z[H_1(X)] \to \F$ be a ring homomorphism with values in a commutative field.
Let $\omega$ be a homological orientation of $X$ and let $\xi \in \Eul(X)$.

Make some intermediate choices: 1) Choose an $o$rdering of the cells
of $X$ and an $o$rientation for each of them. 2) Choose a fundamental
family of cells $\widehat{e}$ which represents $\xi$.
3) Choose a basis $w$ of the $\R$-vector space $H_*(X;\R)$ inducing the orientation $\omega$.
There is the $\F$-complex
$$ 
C_*^\varphi(X):= C_*(\widehat{X})\otimes_{\Z[H_1(X)]} \F
$$
with homology $H^\varphi_*(X) := H_*\left( C_*^\varphi(X)\right)$.
Choices 1 and 2 determine a basis of $C_*(\widehat{X})$ 
with respect to the action of $H_1(X)$ and, so, a basis of 
$C_*^\varphi(X)$ denoted by $\widehat{e}_{oo}$.
Choice 1 also determines a basis of $C_*(X;\R)$ denoted by $oo$. Set
$$
\tau^\varphi(X;\xi,\omega):=
\left\{ \begin{array}{ll}
\sgn\left( \tau(C_*(X;\R); oo, w) \right) \cdot
\tau \left( C_*^\varphi(X); \widehat{e}_{oo}, \varnothing \right)
& \hbox{ if } H_*^\varphi\left(X\right)=0,\\
0 & \hbox{ otherwise.}
\end{array}\right.
$$
It is easily verified that the quantity $\tau^\varphi(X;\xi,\omega)\in \F$ does not depend
on the intermediate choices 1, 2 and 3. Neither it depends on the choice of the base point
$\star \in X$ (which has been introduced to determine the maximal Abelian covering of $X$).
The way $\tau^\varphi(X;\xi,\omega)$ depends on $\xi$ and $\omega$ is given by the rules
$$
\begin{array}{rcll}
\tau^\varphi(X; \xi,-\omega) &=& -\tau^\varphi(X; \xi,\omega)& \\
\tau^\varphi(X; \xi + x,\omega) &=& \varphi(x) \cdot \tau^\varphi(X ; \xi,\omega)
& \ \ \forall x\in H_1(X).
\end{array}
$$

\begin{example}
\label{ex:Milnor--Turaev_torsion}
An important example is provided by the canonical homomorphism
$$
\varphi:\Z[H_1(X)] \to Q\left(\Z[H_1(X)/\Tors\ H_1(X)]\right).
$$
In this case, $\tau^\varphi(X;\xi,\omega)$ is  refered to as the \emph{Milnor--Turaev torsion}
(or, \emph{MT torsion}) of $X$ equipped with $\xi$ and $\omega$. 
Up to multiplication by an element of $\pm H_1(X)/\Tors\ H_1(X)$, 
this fraction happens to coincide with the \emph{Alexander function} of $X$, 
which is defined as the alternated product of the orders 
of the Alexander modules of $X$ if none of them vanishes, and is $0$ otherwise. 
\end{example}

The Abelian group $H_1(X)$ being finitely generated, the quotient ring of 
$\Z[H_1(X)]$ splits (in a unique way) as a product of finitely many fields $\F_i$, see \cite{Turaev_first}. 
The corresponding projections are denoted by $\varphi_i: Q(\Z[H_1(X)]) \to \F_i$.

\begin{definition}
The \emph{Reidemeister--Turaev torsion} (or, \emph{RT torsion}) of the finite connected CW-complex $X$,
equipped with the Euler structure $\xi$ and the homological orientation $\omega$, is
$$
\tau(X;\xi,\omega) := \sum_i \tau^{\varphi_i}(X;\xi,\omega) \
\in \ \bigoplus_i \F_i = Q(\Z[H_1(X)]).
$$
\end{definition}

\subsection{RT torsion of a smooth manifold}

\label{subsec:smooth}

The passing from the category of CW-complexes 
to the category of smooth manifolds can be sketched as follows. 

\subsubsection{Definition}

If $Y$ is a cellular subdivision of a finite connected CW-complex $X$ (such that $\chi(X)=0$),
then there is the \emph{subdivision operator}
$$
\sigma(X,Y): \Eul(X) \to \Eul(Y), \quad [c] \mapsto 
\left[c+ \sum_{\beta,\textrm{ cell of } Y} (-1)^{\dim(\beta)} \cdot \gamma_{\beta}\right] 
$$
where $\gamma_\beta$ is a path contained in the unique open cell $\sigma(\beta)$ of $X$ in which
$\beta$ sits, and connects the center of $\sigma(\beta)$ to that of $\beta$. 
This operator respects the hierarchy of the subdivisions and the RT torsion.

Thus, using triangulations, Turaev proves that the notions of Euler structure and RT torsion extend to polyhedra.
Finally, using smooth triangulations, he extends  those notions from polyhedra to smooth manifolds. 
So, any smooth compact connected $n$-manifold $M$
(such that $\chi(M)=0$) has a $H_1(M)$-affine space of \emph{combinatorial Euler structures}
$$
\Eul_{\hbox{\footnotesize c}}(M)
$$
which, for any smooth triangulation $(X,\rho)$ of $M$, 
can be identified to $\Eul(X)$ via a canonical map  $\rho_*$.
Equipped with a $\xi\in \Eul_{\hbox{\footnotesize c}}(M)$ and a homological orientation $\omega$, 
$M$ has a \emph{RT torsion}
$$
\tau(M;\xi,\omega) \in Q\left(\Z[H_1(M)]\right)
$$
which is equal to $\rho_* \tau(X;\rho_*^{-1}\xi,\rho_*^{-1}\omega)$.

\begin{remark}
\label{rem:cell}
As it is often used, the last two sentences apply verbatim to any cell decomposition $(X,\rho)$ of $M$,
as soon as it can be subdivided to a smooth triangulation.\footnote{
For all the details and proofs that are lacking in this sketchy \S \ref{subsec:smooth}, 
the reader is advised to consult \cite{Turaev_smooth}. 
In particular, Remark \ref{rem:cell} follows from \cite[Lemma 4.2]{Turaev_smooth}.}
\end{remark}

\subsubsection{Geometric Euler structures}

\label{subsubsec:smooth}

Let $M$ be a smooth compact connected $n$-manifold such that $\chi(M)=0$.

A \emph{geometric Euler structure} of $M$ is a non-singular (that is, nowhere zero) vector field on $M$
up to punctured homotopy (that is, up to homotopy among non-singular vector fields 
on $M\setminus \star$, where $\star \in M$). 
By obstruction theory and Poincar\'e duality, the set
$$
\Eul_{\hbox{\footnotesize g}}(M)
$$
of geometric Euler structures is a $H_1(M)$-affine space. Turaev has shown that geometric Euler structures
can be canonically identified with combinatorial ones. 
More precisely, for any smooth triangulation $(X,\rho)$ of $M$, 
there exists a canonical explicit affine map 
$\rho_*:\xymatrix{\Eul(X) \ar[r] &  \Eul_{\hbox{\footnotesize g}}(M)}$ such that, 
if $(X',\rho')$ is another smooth triangulation, 
the following diagram then commutes:
$$
\xymatrix{
\Eul(X) \ar[rd]_-\simeq^-{\rho_*}  \ar[dd]^-\simeq_-{(\rho'^{-1}\circ \rho)_*} & \\
&\Eul_{\footnotesize \rm{g}}(M).\\
\Eul(X') \ar[ru]^-\simeq_-{\rho_*'}  & 
}
$$
(Here, the vertical map is the identification arising from the theory of smooth triangulations,
PL topology and subdivision operators, as alluded to in the previous paragraph.)
So, one gets a canonical affine isomorphism
\begin{equation}
\label{eq:Turaev_identification}
\xymatrix{\Eul_{\hbox{\footnotesize c}}(M) \ar[r]^-\simeq &  \Eul_{\hbox{\footnotesize g}}(M)}.
\end{equation}
In the next sections, geometric Euler structures  will be 
freely identified to combinatorial ones, obtaining thus the set $\Eul(M)$ of \emph{Euler structures} of $M$.

\subsubsection{The case of closed oriented $3$-manifolds}

To come back to our study object, let $M$ be a closed oriented connected $3$-manifold.
The given orientation of $M$ induces an orientation $\omega_M$ of $H_*(M;\R)$:
$$
\omega_M := \left[ \left([\star], b , b^\sharp, [M] \right) \right]
$$
where $[\star] \in H_0(M;\R)$ is the class of a point, $b$ is an
arbitrary basis of  $H_1(M;\R)$,  $b^\sharp$ is the dual basis of $H_2(M;\R)$ 
with respect to the intersection pairing and $[M]$ is the fundamental class.

\begin{definition}
The \emph{RT torsion} of the closed oriented connected $3$-manifold $M$, equipped with the Euler structure $\xi$, is
$$
\tau(M;\xi) := \tau(M; \xi, \omega_M) \in Q\left(\Z[H_1(M)]\right).
$$
\end{definition}

\begin{remark}
\label{rem:Hauptvermutung}
In this definition, we need a priori the manifold $M$ to come with a smooth structure.
The $3$-dimensional Hauptvermutung, together with the fact that two homotopic diffeomorphisms
act the same way on Euler structures, imply that the notions of Euler structure
and RT torsion descend to the topological category.
\end{remark}

\subsection{With the help of Morse theory}

The techniques of smooth triangulations play a  crucial role in the above definition 
of the RT torsion of a smooth manifold $M$ \cite{Turaev_smooth}.
Nevertheless, in practice, it is often much more convenient to work with an arbitrary handle decomposition
given, say, by a Morse function $f:M \to \R$, rather than with a triangulation.
By Remark \ref{rem:cell}, the RT torsion of $M$ can be computed from the cell decomposition 
that is induced by the handle decomposition. But, we still need to understand
how Euler structures relative to such a cell decomposition become geometric through
Turaev's map (\ref{eq:Turaev_identification}). This is the content of Lemma \ref{lem:Morse} below.
Our discussion follows Hutchings--Lee \cite[\S 2]{HL}. See also Burghelea--Haller \cite[\S 3]{BH}.

Let $M$ be a smooth compact connected $n$-manifold such that $\chi(M)=0$.
Let $v$ be a vector field on $M$ with non-degenerate zeros. An \emph{Euler chain} in $M$
\emph{relative} to $v$ is a singular $1$-chain $c$ with boundary 
$$
v^{-1}(0) :=  \sum_{p,\ \hbox{\footnotesize zero of } v} \ind_p(v) \cdot p.
$$
When considered up to homology, Euler chains relative to $v$ form a set 
$$
\Eul(M,v)
$$
which, with the obvious action, is a $H_1(M)$-affine space.

Given two vector fields $v$ and $w$ on $M$ with non-degenerate zeros, 
one can consider the trace of the zeros during a non-degenerate homotopy from $v$ to $w$.
This is a singular $1$-chain in $M$ with boundary $w^{-1}(0) - v^{-1}(0)$,
and it does not depend on the choice of the homotopy up to addition of some boundaries.
This equivalence class is denoted by $w-v$, see \cite[\S 2.2]{BH} for details. 
Then, one gets a canonical affine map $\beta_{v,w}: \Eul(M,v) \to \Eul(M,w)$ 
defined by $[c] \mapsto [c+(w-v)]$.

\begin{lemma}[See \cite{HL,BH}]
\label{lem:alpha}
Let $v$ be a vector field on $M$ with non-degenerate zeros. 
There exists a canonical  affine map $\alpha_v: \Eul(M,v) \to \Eul_{\footnotesize \rm{g}}(M)$
such that, for any other vector field $w$ on $M$ with non-degenerate zeros, the following diagram commutes:
$$
\xymatrix{
\Eul(M,v) \ar[rd]^-{\alpha_v}_-\simeq  \ar[dd]^-\simeq_-{\beta_{v,w}} & \\
&\Eul_{\footnotesize \rm{g}}(M) .\\
\Eul(M,w) \ar[ru]_-{\alpha_w}^-\simeq  & 
}
$$
Moreover, if $c$ is an Euler chain in $M$ relative to $v$ contained in a ball $B\subset M$ and if
$v_c$ is a non-singular vector field on $M$ which coincides with $v$ outside $B$, then $\alpha_v([c])=[v_c]$. 
\end{lemma}
\begin{proof}
Pick a non-singular vector field $v_0$ on $M$. Using the action of $H_1(M)$
on $\Eul_{\footnotesize \rm{g}}(M)$, define $\alpha_{v_0}$ by $[d] \mapsto [v_0] + [d]$.
If $w_0$ is another non-singular vector field on $M$ and if $\alpha_{w_0}$ 
is defined similarly by $[d] \mapsto [w_0] + [d]$, then $\alpha_{w_0}\circ \beta_{v_0,w_0} = \alpha_{v_0}$ since
the homology class $[v_0]-[w_0]$ (given by the affine action of $H_1(M)$ on $\Eul_{\footnotesize \rm{g}}(M)$)
coincides with $w_0-v_0$. Consequently, the map $\alpha_v:= \alpha_{v_0} \circ \beta_{v,v_{0}}$ 
does not depend on the above choice of $v_0$. It has the properties announced by the first statement of the lemma. 

As for the second statement, it suffices to observe that $\beta_{v,v_c}$ sends $[c]$ to $[0]$.
Indeed, one finds a non-degenerate homotopy from $v$ to $v_c$ which is the identity on $M\setminus B$, 
so that $v_c-v$ can be represented by a chain contained in $B$. But, $B$ has trivial homology.
\end{proof}

\begin{remark}
\label{rem:link_with_Turaev}
Let $(X,\rho)$ be a smooth triangulation of $M$ and let $v(X,\rho)$ 
be the Whitney singular vector field on $M$ associated to this triangulation \cite{HT}.
Then, the composition
$$
\xymatrix{
\Eul(X) \ar@{=}[r]^-{\rho} & \Eul\left(M,v(X,\rho)\right) \ar[rr]_-\simeq^-{\alpha_{v(X,\rho)}} & &
 \Eul_{\footnotesize \rm{g}}(M)
}
$$
is Turaev's map $\rho_*$ that has been mentioned in \S \ref{subsubsec:smooth}.
\end{remark}

\begin{lemma}[Hutchings--Lee \cite{HL}] 
\label{lem:Morse}
Let $f:M\to \R$ be a Morse function together
with a Riemannian metric, such that the Smale condition is satisfied.
Let $\left(X_{f},\rho_{f}\right)$ be the Thom--Smale cell decomposition of $M$ 
associated\footnote{Recall that the Smale condition requires that, for any critical points
$x$ and $y$ of $f$, the ascending manifold of $x$ is transverse to the descending manifold
of $y$. The open $i$-cells of the Thom--Smale cell decomposition of $M$
are the descending manifolds from index $i$ critical points. See \cite{Laudenbach}
for a precise description of this CW-complex.} 
to $f$, and let $\nabla f$ be the gradient field of $f$ with respect to the given metric. 
Then, the diagram
$$
\xymatrix{
\Eul\left(X_f\right) \ar@{=}[r]^-{\rho_f} \ar[d]^-\simeq_-{\rho_{f,*}}
& \Eul(M,\nabla f) \ar[d]_-\simeq^-{\alpha_{\nabla f}} \\ 
\Eul_{\footnotesize \rm{c}}(M) \ar[r]^\simeq &  \Eul_{\footnotesize \rm{g}}(M),
}
$$
where the bottom map is Turaev's identification (\ref{eq:Turaev_identification})
and where the map $\rho_{f,*}$ is given by Remark \ref{rem:cell}, commutes.
\end{lemma}

\begin{proof} This is sketched in the proof of \cite[Lemma 2.13]{HL}, whose argument 
we would like to develop.

\begin{claim}
\label{claim:Morse}
Let $f_1,f_2: M\to \R$ be Morse functions, each  with a Riemannian metric 
such that the Smale condition is satisfied. Then, the following diagram commutes: 
\begin{equation}
\label{diag:Morse}
\xymatrix{
\Eul\left(X_{f_1}\right) \ar@{=}[r]^-{\rho_{f_1}} 
\ar[d]^-\simeq_-{ \rho_{f_2,*}^{-1} \circ \rho_{f_1,*} }
& \Eul\left(M,\nabla f_1\right) \ar[d]_-\simeq^-{\beta_{\nabla f_1,\nabla f_2}}\\
\Eul\left(X_{f_2}\right) \ar@{=}[r]_-{\rho_{f_2}} & \Eul\left(M,\nabla f_2\right).
}
\end{equation}
\end{claim}

\noindent
Assuming this, let $(X,\rho)$ be an arbitrary smooth triangulation of $M$. 
There exists a Morse function $f_1: M \to \R$ \emph{compatible} with $(X,\rho)$,
in the sense that $\nabla f_1$ coincides essentially 
with the Whitney singular vector field $v(X,\rho)$.
Then, the two cell decompositions $(X_{f_1}, \rho_{f_1})$ and $(X,\rho)$ coincide. 
See Figure \ref{fig:PL_to_Morse}.
\begin{figure}[h]
\centerline{\relabelbox \small 
\epsfxsize 3truein \epsfbox{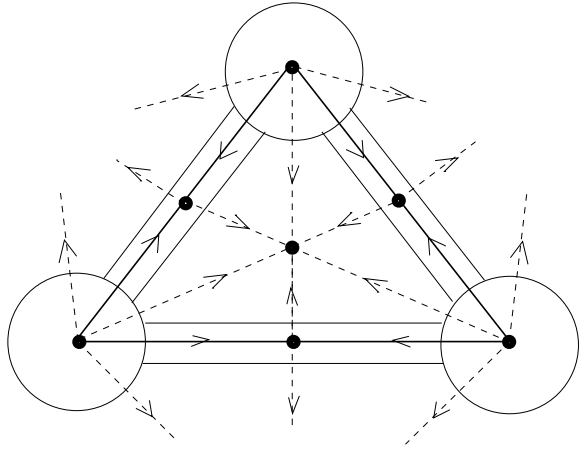}
\endrelabelbox}
\caption{The triangulation $(X,\rho)$, the critical points of $v(X,\rho)$
with some of its flow lines, and the handle decomposition induced
 by a compatible Morse function  ($n=2$).}
\label{fig:PL_to_Morse}
\end{figure}
Applying the above claim to this $f_1$ and to $f_2:=f$, 
we get the internal square of the following commutative diagram:
$$
\xymatrix{
& \Eul(X) \ar[dd] \ar@{=}[r]^-\rho \ar[ld]_-{\rho_*} &
 \Eul\left(M,v(X,\rho)\right) \ar[dd]^-{\beta_{v(X,\rho),\nabla f}} \ar[rrd]^-{\alpha_{v(X,\rho)}} && \\
\Eul_{\footnotesize \rm{c}}(M)  &&&& \Eul_{\footnotesize \rm{g}}(M) \\
&\Eul(X_f) \ar@{=}[r]_-{\rho_f} \ar[lu]^-{\rho_{f,*}} & \Eul(M,\nabla f) \ar[rru]_-{\alpha_{\nabla f}}&&
}
$$
The triangle right is given by Lemma \ref{lem:alpha}. 
The conclusion follows from Remark \ref{rem:link_with_Turaev}.\\

We now prove Claim  \ref{claim:Morse} following Laudenbach's bifurcation analysis\footnote{This analysis 
requires the metric to have a special form near the critical points
of the Morse function. It is enough to prove the claim assuming this special Morse condition for $f_1$
and $f_2$. Nevertheless, we simplify the exposition by not taking care of the metric.} 
of the Thom--Smale complex \cite{Laudenbach}.
By Cerf theory, there exists a path $(f_t)_{t}$ of pairs 
(a smooth function $M\to \R$, a metric) connecting $f_1$ to $f_2$, and such that 
$f_t$ is a Morse function satisfying the Smale condition at each time $t$, 
except in an $\varepsilon$-neighborhood of a finite number of times $1 < t_1< \cdots < t_r < 2$. 
Around each time $t_k$, 
one of the following scenari may occur:
\begin{itemize}
\item[(a)] From $f_{t_k-\varepsilon}$ to $f_{t_k+\varepsilon}$, there is the birth or the death
of two critical points of consecutive indices $i,i+1$.
\item[(b)] From $f_{t_k-\varepsilon}$ to $f_{t_k+\varepsilon}$, 
the functions satisfy the Morse condition, 
there is a regular level $L:= f^{-1}_{t_k}(a)$ below which and above which the Smale condition is satisfied
for the two cobordisms it delimitates. But, at time $t_k$ and at the level $L$, 
the trace of an ascending manifold from a critical point $x$ of index $i$ (in the lower cobordism) 
fails to be tranverse in a unique point $p$ to the trace of the descending manifold of a critical point $y$ of index $j$ 
(in the upper cobordism). (This failure can be of two types: $j>i$ or  $j=i$.)
\end{itemize}
In terms of handle decompositions of $M$, (a) is a stabilization or a destabilization, and
(b) results from an isotopy of the attaching region of an index $j$ handle, 
when its lower sphere crosses the upper sphere of an index $i\leq j$ handle.
(This is a handle sliding, when $i=j$). 

If $r=0$, there exists an ambiant isotopy $(\phi_t)_t$ from $\phi_1=\Id_M$ to a certain $\phi_2$
such that $f_t=\phi_t\circ f_1$ for each time $t$.
Let $d$ be the trace of the critical points of $f_1$ during
the isotopy. For any $[c]\in \Eul(X_{f_1})$, 
the combinatorial Euler structures of $M$ $\rho_{f_1,*}([c])$ and $\rho_{f_2,*}\left([\phi_2(c)]\right)$ are equal; 
but, $\phi_2(c)$ is homologous to $c+ d$ (performing the isotopy $(\phi_t(c))_t$ backwards). 
Besides, the homotopy $(f_t)_t$ between  $f_1$ and $f_2$ 
induces a non-degenerate homotopy $(\nabla f_t)_t$ from $\nabla {f_1}$ to $\nabla {f_2}$,
so $d$ represents $\nabla f_2 - \nabla f_1$. We conclude that the diagram (\ref{diag:Morse}) commutes.

So, it suffices to consider the case when $r=1$ and  $f_1\leadsto f_2$ 
by a one-event scenario of type (a) or (b). 
In the (b) case, the same isotopy argument as above applies. In the (a) case,
we assume that a stabilization occurs $\varnothing \leadsto A^i \cup A^{i+1}$
with the birth of two new handles. After having possibly isotoped handles of indices $\geq i$,
we can assume that it happens in the interior of a $n$-handle $A^n$. See Figure \ref{fig:stabilization}.
\begin{figure}[h]
\centerline{\relabelbox \small 
\epsfxsize 3.5truein \epsfbox{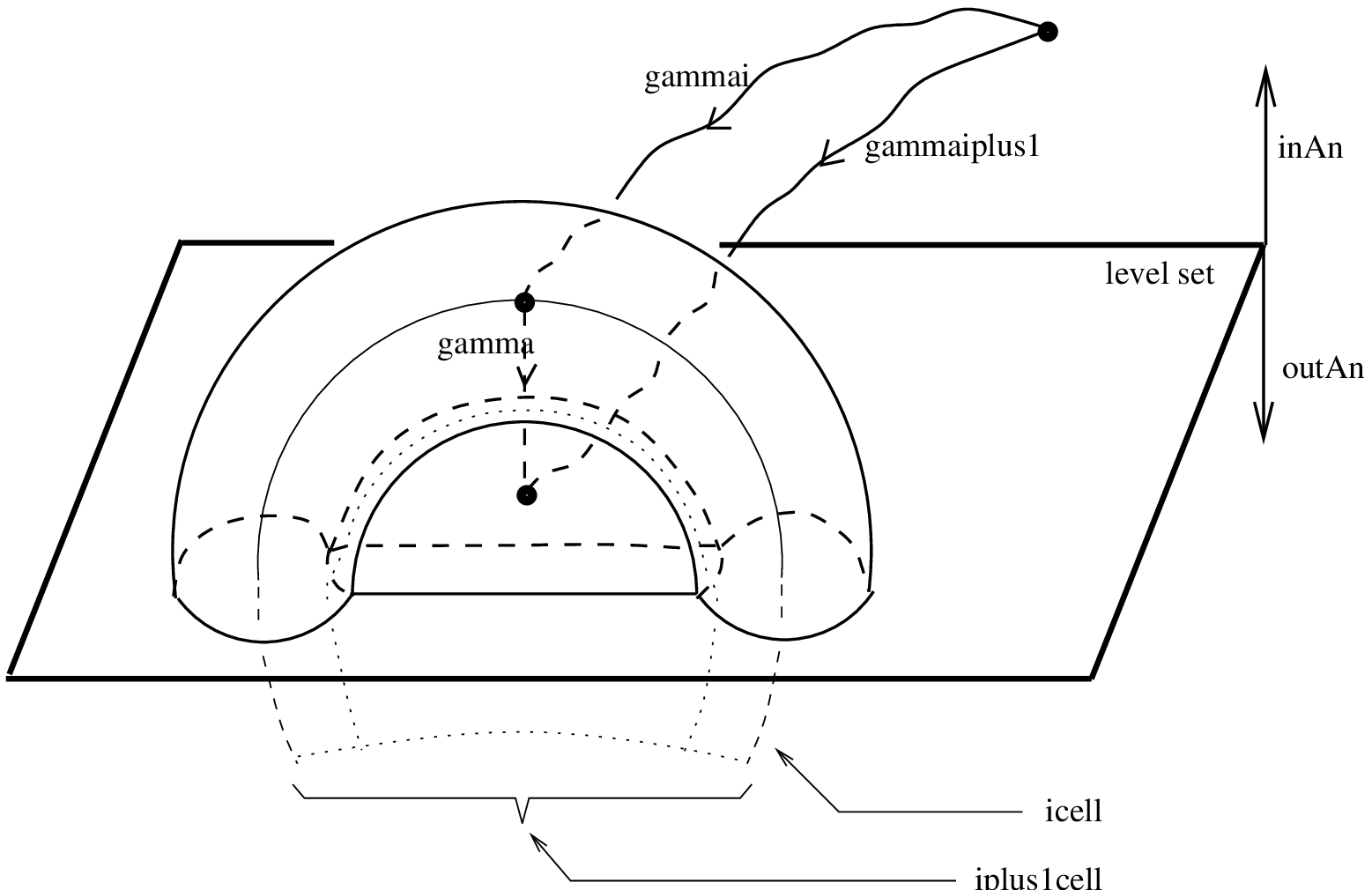}
\adjustrelabel <0.3cm,-0.1cm> {gamma}{$\gamma$}
\adjustrelabel <0.1cm,-0.1cm> {gammai}{$\gamma_i$}
\adjustrelabel <-0.1cm,-0.2cm> {gammaiplus1}{$\gamma_{i+1}$}
\adjustrelabel <0.1cm,0cm> {inAn}{in $A^n$}
\adjustrelabel <0.cm,0cm> {outAn}{out $A^n$}
\adjustrelabel <0.cm,0cm> {icell}{the extra $i$-cell}
\adjustrelabel <0.cm,0cm> {iplus1cell}{the extra $(i+1)$-cell}
\adjustrelabel <-0.7cm,-0.1cm> {level set}{a level set} 
\endrelabelbox}
\caption{A stabilization is a subdivision ($n=3, i=1$).}
\label{fig:stabilization}
\end{figure}
Then, $(X_{f_2}, \rho_{f_2})$ is a subdivision of $(X_{f_1}, \rho_{f_1})$ 
with one extra $i$-cell  and one extra $(i+1)$-cell corresponding to the 
descending manifolds associated to $A^i$ and $A^{i+1}$ respectively.
They subdivide the $n$-cell emerging from $A^n$.
Let $\gamma_i$ and $\gamma_{i+1}$ be paths in $A^n$ connecting its center to the centers 
of $A^i$ and $A^{i+1}$ respectively:  
The subdivision operator $\sigma\left(X_{f_1},X_{f_2}\right)$ sends 
$[c]$ to $[c+ (-1)^i \gamma_i +(-1)^{i+1} \gamma_{i+1}]$. 
On the other hand, the homotopy $(f_t)_t$ induces a non-degenerate homotopy from 
$\nabla f_1$ to $\nabla f_2$ which is fixed outside a neigborhood of $A^i \cup A^{i+1}$. 
The trace of the critical points during this homotopy is a small arc $\gamma$ 
connecting the center of $A_i$ to the center of $A_{i+1}$ (with the convenient orientation, 
according to the parity of $i$). 
Since $(-1)^i \gamma_i+ (-1)^{i+1}\gamma_{i+1}$ is homologous to $\gamma$ in the ball $A^n$,
the commutativity of (\ref{diag:Morse}) follows.
\end{proof}

\section{Finite-type invariants of manifolds with structures}
\label{sec:FTI_with_structures}

In this section, we settle the context in which finiteness properties
of the RT torsion will be proved. 

\subsection{Canonical correspondences induced by a Torelli surgery}

\label{subsec:Torelli}

Let $M$ be a closed oriented connected $3$-manifold. Given a handlebody $H \subset M$
and a Torelli automorphism $h$ of $\partial H$, we define
$$
M_h := (M \setminus \int\ H) \cup_h H
$$
and call the move $M \leadsto M_h$ a \emph{Torelli surgery}.

The move $M \leadsto M_h$ induces a canonical 
isomorphism in homology. This is the only map $\Phi_h$ which makes the diagram
\begin{equation}
\label{diag:Phi}
\xymatrix{
&H_1(M)\ar@{->}[dd]^{\Phi_h}_\simeq\\
H_1\left(M\setminus \int\ H \right) \ar@{->>}[ru]^{\hbox{\footnotesize incl}_*}
\ar@{->>}[rd]_{\hbox{\footnotesize incl}_*} & \\
&H_1\left(M_h\right)
}
\end{equation}
commute. (It exists and is injective because $h$ acts trivially in homology,
so that a homology class 
in $M\setminus \int\ H$ vanishes in $M$ if and only if it does in $M_h$.)

Dually, the move $M \leadsto M_h$ induces a canonical
correspondence between parallellizations or, equivalently, Spin-structures.
A \emph{parallelization} of $M$ is a trivialization of $\T M$, the oriented tangent bundle of $M$,
up to punctured homotopy. There is a unique map $\Theta_h$ which makes the diagram
\begin{equation}
\label{diag:Theta}
\xymatrix{
\Parall(M) \ar[dd]_{\Theta_h}^\simeq \ \ \ \ar@{>->}[rd]^{\hbox{\footnotesize incl}^*} &\\
& \Parall\left(M\setminus \int\ H\right)\\
\Parall(M_h) \ \ \ \ \ar@{>->}[ru]_{\hbox{\footnotesize incl}^*} &
}
\end{equation}
commute \cite{Massuyeau}. (It exists and is surjective because $h$ 
acts trivially at the level of Spin-structures -- as follows from \cite{Johnson_spin} --
so that a parallelization of $M\setminus \int\ H$
extends to $M$ if and only if it does to $M_h$.)

Furthermore, the Torelli surgery $M \leadsto M_h$ induces a canonical
correspondence between Euler structures. To see this, we need the canonical
map $\beta:\Parall(M) \to \Eul(M)$ that forgets the second and third vectors 
of a trivialization of $\T M$.  

\begin{lemma}[See \cite{DM}]
\label{lem:Omega} There exists a unique bijection 
$$
\Omega_h: \Eul(M) \longrightarrow \Eul\left(M_h\right), \quad \xi \longmapsto \xi_h,
$$
which is affine over $\Phi_h$ and makes the following diagram commute: 
$$
\xymatrix{
{\Parall}(M)\ar[d]_\beta \ar[r]^-{\Theta_h}_-{\simeq} 
& {\Parall}\left(M_h\right) \ar[d]^\beta\\
\Eul(M) \ar[r]^-{\simeq}_-{\Omega_h} & \Eul(M_h).
}
$$
\end{lemma}

\begin{proof}
Pick a parallelization $\sigma_0$ of $M$. Then, $\Omega_h$ must be defined by 
$\Omega_h(\beta(\sigma_0)+x) := \beta \Theta_h(\sigma_0)+\Phi_h(x)$
for any $x\in H_1(M)$. We must verify that the choice of $\sigma_0$ does not matter.
Recall that, by obstruction theory, the set $\Parall(M)$ is an affine space over $H^1(M;\Z_2)$.
It follows from the above definitions that $\Theta_h$ is affine
over the inverse of $\Phi_h^*: H^1(M_h;\Z_2) \to H^1(M;\Z_2)$.
Also, the map $\beta$ is affine over the Bockstein homomorphism $\beta$.
So, that verification amounts to check that the diagram
$$
\xymatrix{
H^1(M;\Z_2) \ar[d]_-\beta & H^1(M_h;\Z_2) \ar[l]_-{\Phi_h^*} \ar[d]^-\beta\\
H^2(M) \ar[r]_-{P \Phi_h P^{-1}}& H^2(M_h)
}
$$
commutes. We expand it as follows:
$$
\xymatrix{
H^1(M;\Z_2) \ar[ddd]_-\beta &&& 
H^1(M_h;\Z_2) \ar[lll]_-{\Phi_h^*} \ar[ddd]^-\beta\\
&H_2(M;\Z_2) \ar[d]_-\beta \ar[lu]^-{P} & H_2(M_h;\Z_2) \ar[d]^-\beta \ar[ru]_-{P} &\\
&H_1(M) \ar[ld]_-{P} \ar[r]^-{\Phi_h} & H_1(M_h) \ar[rd]^-P &\\
H^2(M) \ar[rrr]_-{P \Phi_h P^{-1}} &&& H^2(M_h).
}
$$
One proves that the above hexagon commutes using the fact 
that the isomorphism $\Phi_h$ preserves the linking pairings.
\end{proof}

\begin{remark}
\label{rem:smoothing}
Technically, in what precedes, we had to fix a smooth structure on $M$. Next,
we chose on $M_h$  a smooth structure which induces on $M \setminus \int\ H$  
and $H$ those ones induced by $M$. 
(Such a smooth structure on $M_h$ is not unique, but it is up to a diffeomorphism 
which is the identity on $M \setminus \int\ H$ and which is homotopic to the identity.) 
One easily checks that, in the sense of Remark \ref{rem:Hauptvermutung},
the map $\Omega_h$ is independent of those two consecutive choices of smooth structures.
\end{remark}

\begin{remark}
\label{rem:canonicity}
Let $f$ be a diffeomorphism of $\partial H$ which extends to $H$ and acts trivially in homology.
Let $\phi: M_h \to M_{hf}$ be a diffeomorphism which is the identity on 
$M\setminus \int\ H$. Diagrams (\ref{diag:Phi})
and (\ref{diag:Theta}) imply that $\phi_* \Phi_h= \Phi_{hf}$ and 
$\phi_* \Theta_h= \Theta_{hf}$ respectively. So, we have that
$\phi_* \Omega_h= \Omega_{hf}$ as well.

Similarly, let $(\phi_t)_{t\in [0,1]}$ be an ambiant isotopy of $M$ sending $H$ to $H':=\phi_1(H)$.
Then, $h':= \phi_1| \circ h \circ \phi_1|^{-1}$ is a Torelli automorphism of $\partial H'$
and $\phi_1$ induces a diffeomorphism $\phi: M_h \to M_{h'}$.
Using the same diagrams, one checks successively that $\phi_* \Phi_h= \Phi_{h'}$,
$\phi_* \Theta_h= \Theta_{h'}$ and $\phi_* \Omega_h= \Omega_{h'}$.
\end{remark}

The correspondence $\Omega_h$ has been defined equivalently in \cite[\S 3.2.1]{DM}
using some kinds of relative Euler structures and their gluings.
Next lemma gives a third description of the bijection $\Omega_h$ 
that makes gluing of vector fields quite explicit. It needs the \emph{Chil\-ling\-worth homomorphism}
$$
t: \mathcal{T}(\Sigma_{g,1}) \longrightarrow 2\cdot H_1(\Sigma_{g,1}) \subset H_1(\Sigma_{g,1}), 
\quad f \longmapsto t(f)
$$
defined on the Torelli group of a compact connected oriented surface $\Sigma_{g,1}$ 
of genus $g$  with $1$ boundary component. The homology class $t(f)$ corresponds to the obstruction
$$
\dd f^{-1}(s) - s \in H^1(\Sigma_{g,1})\simeq H_1(\Sigma_{g,1},\partial \Sigma_{g,1}) \simeq H_1(\Sigma_{g,1})
$$ 
to homotope a non-singular vector field $s$ on $\Sigma_{g,1}$ (any one) to its image  under $f^{-1}$.
It happens to be even. See \cite[\S 5]{Johnson_homomorphism} for details.

\begin{remark}
\label{rem:Chillingworth}
The Chillingworth homomorphism is easily computable, 
for instance using Fox's free differential calculus. More precisely, let
$(z_1,z_2,\dots,z_{2g-1},z_{2g})$ be a symplectic basis of 
$\pi:=\pi_1\left(\Sigma_{g,1},\star\right)$ where $\star \in \partial \Sigma_{g,1}$.
The \emph{Magnus representation} of the Torelli group of $\Sigma_{g,1}$ is the group homomorphism
$$
\mathcal{T} \xymatrix{\left(\Sigma_{g,1}\right) \ar[r]^-{r^{\mathfrak{\scriptsize a}}} & 
\hbox{GL}\left(2g;\Z[\pi/\pi']\right)
}
$$
defined by 
$$
r^{\mathfrak{\scriptsize a}}(f) =\mathfrak{a}
\left(\overline{\frac{\partial f_*(z_j)}{\partial z_i}}\right)_{i,j}.
$$
Here, $\mathfrak{a}$ denotes the ring homomorphism induced by the Abelianization 
$\pi \to \pi/\pi'$ while the bar denotes the ring anti-homomorphism induced by the inversion of $\pi$
\cite[\S 5]{Morita}.  Combining results from Johnson \cite[Theorem 2]{Johnson_homomorphism}
and Morita \cite[Theorem 6.1]{Morita}, one gets the following formula: 
$$
\begin{array}{rcll}
t(f) &=& 2g - \Trace\left(r^{\mathfrak{a}}(f)\ \mod\ I(\Z[\pi/\pi'])^2\right)&\\
&=& \sum_{i=1}^{2g} \left(\mathfrak{a}\left(\frac{\partial f_*(z_i)}{\partial z_i}\right)-1\right) 
&\in I(\Z[\pi/\pi']) / I(\Z[\pi/\pi'])^2 \simeq \pi/\pi'.
\end{array}
$$
\end{remark}

\begin{lemma}
\label{lem:Omega_more_concretely}
Assume, after an isotopy, that $h$ is the identity on a disk $D\subset \partial H$.
Suppose given a bi-collar neighborhood $[-1,1] \times \partial H \hookrightarrow M$
of $\partial H$ in $M$. Choose a smooth structure on $M$, and endow $M_h$ 
with the unique smooth structure compatible with those on 
$M\setminus \int\ H$ and $H$ induced by $M$ and such that, 
for each $x \in \partial H$, the intervals $[-1,0]\times x$ and $[0,1] \times h(x)$ piece
together to give a smooth interval in $M_h$.

Any Euler structure $\xi$ of $M$ can be represented by a non-singular vector field $v$
which is outwards (respectively inwards) normal to $H$ on $\partial H\setminus D$ --
i.e. coincides with $\varepsilon \cdot \partial/\partial t$ on $0\times (\partial H\setminus D)$
where $\varepsilon=+1$ (respectively $\varepsilon=-1$).
Then, for such a vector field $v$, the non-singular vector field 
$$
v_h:= v|_{M\setminus \int H} \cup_h v|_H
$$
of $M_h$ represents the Euler structure 
$$
\xi_h+ \varepsilon\cdot  \incl_*\left(t(h|_{\partial H \setminus D})/2\right).
$$
\end{lemma}

\begin{proof}
Using the fact that the correspondence $\Omega_h$ sends an Euler structure that
comes from a parallelization to an Euler structure with the same property,
one easily checks that $\Omega_h$ commutes
with the involution of Euler structures defined by $[u]\mapsto [-u]$. Consequently,
it is enough to prove the lemma in the outwards case $\varepsilon=+1$.

Let $\xi \in \Eul(M)$. We first prove that there exists a representant $v$ of $\xi$ 
which is outwards normal to $H$ on $\partial H\setminus D$.
In general, for $N$ a compact oriented $3$-manifold with boundary and
a non-singular section $s$ of T$N|_{\partial N}$, let us call
an \emph{Euler structure} on $N$ \emph{relative} to $s$
a punctured homotopy class rel $\partial N$ of non-singular vector fields on $N$ 
that extend $s$ \cite[\S 1.3.3]{DM}.
We denote by $\Eul(N,s)$ the set of such structures. 
Obstruction theory tells us that there is an obstruction 
$w(N,s) \in H^3(N,\partial N)$ to the existence of such structures
and, when the latter happens to vanish, that the set $\Eul(N,s)$ is naturally
a $H_1(N)$-affine space (using Poincar\'e duality). Furthermore, 
obstruction calculi on the double $N\cup_{\rm{\scriptsize Id}} (-N)$
and an application of the Poincar\'e--Hopf theorem lead to
\begin{equation}
\label{eq:obstruction}
2\cdot \langle w(N,s), [N,\partial N]\rangle
= \left\langle e\left(\hbox{T}N|_{\partial N}/ \langle s \rangle\right), [\partial N]\right\rangle \in \Z
\end{equation}
where $e$ denotes the Euler class.
 
Let $(s,s',s'')$ be a trivialization of T$M|_{\partial H}$ such that $s$ is 
the outwards normal vector field to $H$ on $\partial H\setminus D$. 
This exists, as illustrated on Figure \ref{fig:vector_field}.
\begin{figure}[h]
\centerline{\relabelbox \small 
\epsfxsize 4truein \epsfbox{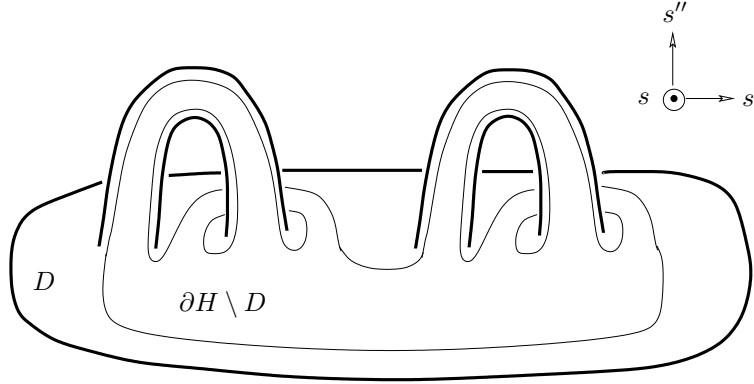}
\adjustrelabel <0.cm,0cm> {D}{$D$}
\adjustrelabel <0.cm,0cm> {C}{$\partial H \setminus D$}
\adjustrelabel <-0.1cm,0cm> {w}{$s$}
\adjustrelabel <-0.1cm,0cm> {w'}{$s'$}
\adjustrelabel <-0.1cm,0cm> {w''}{$s''$}
\endrelabelbox}
\caption{The handlebody $H$ and a trivialization $(s,s',s'')$
of T$H|_{\partial H}$ such that $s$ is the outwards normal vector field on $\partial H\setminus D$.}
\label{fig:vector_field}
\end{figure}
According to (\ref{eq:obstruction}), the obstruction to extend $s$ to $H$ vanishes (which is also apparent
in Figure \ref{fig:vector_field}), and the obstruction to extend $s$ to $M\setminus \int\ H$ does too. 
So, we can consider the obvious gluing map
$$
\xymatrix{
\Eul(H,s) \times \Eul(M\setminus \int\ H, s) \ar[r]^-{\cup} & \Eul(M).
}
$$
This map is affine over 
$\incl_*\oplus \incl_*: H_1(H) \times H_1(M\setminus \int\ H) \to H_1(M)$
and, so, is surjective. This implies the first claim of the lemma.

Next, assume that the second claim of the lemma holds for another Euler
structure $\xi_0$ on $M$ rather than our $\xi$.
Represent $x:=\xi_0-\xi\in H_1(M)$ by a knot $K$ in $M\setminus \int\ H$, 
disjoint from the collar neighborhood of $\partial H$, 
and modify the given $v$ by Reeb surgery along $K$. 
The resulting vector field $v_0$ represents $\xi_0$ so that, by assumption,
$(v_0)_h$ represents $(\xi_0)_h+ \incl_*( t({h|})/2 )$. Since $v_h$ is obtained from
$(v_0)_h$ by Reeb surgery along $-K$, we deduce that it represents
$$
\left((\xi_0)_h+ \incl_*(t({h|})/2 )\right) + [-K] = 
\left((\xi_0)_h - \Phi_h(x)\right) +  \incl_*(t({h|})/2 ) = \xi_h +  \incl_*(t({h|})/2 ),
$$
which shows that the lemma holds for $\xi$ as well. 

Thus, it is enough to prove the second claim of the lemma for a particular $\xi$, 
e.g. one which comes from a parallelization $\sigma$.
Then, there exist some non-singular vector fields $v'$ and $v''$ on $M$ such that
$(v,v',v'')$ is a trivialization of T$M$ with $v'$ and $v''$ tangent to $\partial H$ on $\partial H \setminus D$.
Since the Torelli automorphism $h^{-1}$ fixes $\sigma|_{\partial H}$ (as follows from \cite{Johnson_spin}), 
there exists a homotopy $F=(f_t)_{t\in [-1,0]}$ from the trivialization $(v,v',v'')$ 
of T$H|_{\partial H}= \R \oplus \hbox{T} \partial H$ to its image
$$ 
(\Id \oplus \dd h^{-1}) (v,v',v'')=
\left(v,(\Id \oplus \dd h^{-1})(v') ,(\Id \oplus \dd h^{-1})(v'')\right).
$$
According to the diagram (\ref{diag:Theta}), $\Theta_h(\sigma)$ is represented by
\begin{equation}
\label{eq:parallelization}
(v,v',v'')|_{M\setminus \int\ H}\cup_h
\left(F \cup (v,v',v'')|_{H} \right)
\end{equation}
where the second gluing ``$\cup$" is given by the collar neigborhood of $\partial H$ in $H$.
According to Lemma \ref{lem:Omega},
the Euler structure $\xi_h$ equals $\beta\Theta_h(\sigma)$ and, so, 
is represented by the first vector field of the triplet (\ref{eq:parallelization}). 
This is \emph{not} necessarily $v_h$ because the homotopy $F$ may have moved $v|_{\partial H}$
to itself between times $t=-1$ and $t=0$. 
So, we are interested in the space of relative Euler structures
$$
\Eul\left( [-1,0]\times \partial H, \left((-1) \times v|_{\partial H}\right) \cup \left( 0\times v|_{\partial H}\right) \right)
$$
and its two elements ``$v$'' (the vector field equal to $t\times v|_{\partial H}$ at each time $t\in [-1,0]$) 
and $\beta(F)$. To conclude, we need to compute their difference
$$
\Delta := \textrm{``$v$''} - \beta(F) \in H^2([-1,0] \times \partial H, \partial [-1,0] \times \partial H)
\simeq H_1(\partial H),
$$
since it is such that $\xi_h + \incl_* \circ P^{-1}(\Delta)$ is represented by $v_h$. 
For this, consider the map
$$
\xymatrix{
\Eul\left( [-1,0]\times \partial H, \left((-1) \times v|_{\partial H}\right) 
\cup \left(0\times v|_{\partial H}\right) \right)
\ar[r]^-c & H^2([-1,0] \times \partial H, \partial [-1,0] \times \partial H)
}
$$
that assigns to $[r]$ the obstruction
$$
e\left(\hbox{T}([-1,0]\times \partial H)/ \langle r\rangle,
\left((-1)\times v'|_{\partial H}\right) \cup \left( 0 \times (\Id \oplus \dd h^{-1})(v'|_{\partial H})\right) \right)
$$
to extend the section of the normal bundle of $r$ in $\hbox{T}([-1,0]\times \partial H)$
given by $(-1)\times v'|_{\partial H}$ and $0\times (\Id \oplus \dd h^{-1})(v'|_{\partial H})$ 
on the boundary. The map $c$ is affine over the doubling map. On the one hand, it vanishes on $\beta(F)$.
On the other hand, note that the bundle $\hbox{T}([-1,0]\times \partial H)/ \langle r\rangle$ for $r=\hbox{``$v$''}$ 
is the pull back of $(\R \oplus \hbox{T}\partial H)/ \langle v|_{\partial H} \rangle$ 
by the cartesian projection $[-1,0] \times \partial H \to \partial H$; so, $c(``v\hbox{''})$ 
is the image  by the suspension isomorphism
$$
S: H^1(\partial H) \to H^2([-1,0] \times \partial H, \partial [-1,0] \times \partial H)
$$
of the obstruction to homotope  $v'|_{\partial H}$ to $(\Id \oplus \dd h^{-1})(v'|_{\partial H})$,
which are both sections of $(\R \oplus \hbox{T}\partial H)/ \langle v|_{\partial H} \rangle$.
The latter obstruction corresponds by the transfer map
$$
H^1( \partial H\setminus D, \partial( \partial H\setminus D)) \simeq 
\xymatrix{H_1(\partial H\setminus D) \ar[r]^-{\incl_*}_-\simeq& H_1(\partial H)} \simeq H^1(\partial H)
$$
to the obstruction
$\dd h^{-1}(v')-v' \in H^1( \partial H\setminus D, \partial( \partial H\setminus D))$
to homotope, relatively to the boundary, 
the tangent vector field $v'$ on the surface $\partial H \setminus D$ to $\dd h^{-1}(v')$. 
This is essentially the Chillingworth class of $h|_{\partial H\setminus D}$. 
Since $S$ corresponds to $\incl_*: H_1(\partial H) \to H_1([-1,0] \times \partial H)$
through Poincar\'e dualities, we obtain that
$$
c(``v\hbox{''})= S \circ \underbrace{P \circ \incl_* \circ P^{-1}}_{\hbox{transfer}} \left( \dd h^{-1}(v')-v'\right)
= S \circ P \circ \incl_* ( t(h|) ) = P \circ \incl_* \left(t(h|) \right).
$$
We conclude that 
$$
2\cdot \Delta= 2 \cdot \left(\textrm{``$v$''} - \beta(F) \right)= 
c(\textrm{``$v$''}) - c(\beta(F)) = P\circ  \incl_*\left( t(h|) \right)
$$
which implies that $P^{-1}(\Delta)= \incl_* \left(t(h|)/2\right) \in H_1(\partial H)$.
\end{proof}

\subsection{Adding structures to the Goussarov--Habiro theory}

\label{subsec:GH_with_structures}

Let $\mathcal{M}$ be the set of closed oriented connected $3$-manifolds, 
up to orientation-preserving homeomorphisms. 
Usually, the Goussarov--Habiro theory applies to the set $\mathcal{M}$ \cite{GGP,Habiro,Garoufalidis}.
In this paper, we place it in a more general context where manifolds come with additional structures.

More precisely, we fix a finitely generated Abelian group $G$, and we consider triples of the form
$$
(M,\xi,\psi)
$$
where $M$ is a closed connected oriented $3$-manifold, $\xi$ is an Euler structure on $M$ and 
$\psi: G \to H_1(M)$ is an isomorphism. 
An \emph{isomorphism} between two such objects  $(M_1,\xi_1,\psi_1)$ and $(M_2,\xi_2,\psi_2)$ 
is an orientation-preserving homeomorphism $f:M_1 \to M_2$ which carry $\xi_1$ to $\xi_2$
and such that $f_* \circ \psi_1= \psi_2$. The set of isomorphism classes is denoted by 
$$
\mathcal{ME}(G).
$$
When either Euler structures or homological parametrizations are not taken into account,
we denote the corresponding sets by $\mathcal{M}(G)$ and $\mathcal{ME}$ respectively. 

Because of the canonical correspondences that it induces
(\S \ref{subsec:Torelli}), the Torelli surgery  extends from $\mathcal{M}$ to $\mathcal{ME}(G)$.
Given a triple $(M,\xi,\psi)$ as above, a handlebody $H \subset M$ and a Torelli automorphism $h$ of $\partial H$, 
we call the move $(M,\xi,\psi) \leadsto (M,\xi,\psi)_h$, where  
$(M,\xi,\psi)_h:= (M_h, \Omega_h(\xi), \Phi_h \circ \psi)$, a \emph{Torelli surgery}.

\begin{definition}
\label{def:fti_structures}
Let $A$ be an Abelian group and
let $f:\mathcal{ME}(G) \to A$ be an invariant of closed connected oriented $3$-manifolds with Euler structure
and homology parametrized by $G$. 
This is a \emph{finite-type invariant} of \emph{degree} at most $d$ if,
for any $(M,\xi, \psi) \in \mathcal{ME}(G)$ and for any family $\Gamma$ of
$d+1$ Torelli automorphisms of the boundaries of pairwise disjoint handlebodies in $M$,
the following identity holds:
\begin{equation}
\label{eq:fti_structures}
\sum_{\Gamma' \subset \Gamma} (-1)^{|\Gamma'|} \cdot f\left((M,\xi,\psi)_{\Gamma'}\right) =0 \in A.
\end{equation}
Here, $M_{\Gamma'}$ denotes the manifold obtained from $M$ 
by the simultaneous surgery defined by those elements of  $\Gamma$ shortlisted in $\Gamma'$.
\end{definition}

In order to connect this extension of the Goussarov--Habiro theory with the existing literature,
we now make a few comments which are dedicated to the initiated reader. 

\begin{remark}
The Goussarov--Habiro theory for $\mathcal{M}$ contains non-trivial degree $0$ invariants. 
According to Matveev \cite{Matveev},  they classify the pair
(homology, linking pairing) of a closed connected oriented $3$-manifold. 
This result has been extended to $\mathcal{ME}$ in \cite{DM}. 
Adding the homological parametrizations, one obtains\footnote{Combining Theorem 2 and Remark 3.8 from \cite{DM}.} 
that the ``universal'' degree $0$ invariant of the Goussarov--Habiro theory for 
$\mathcal{ME}(G)$ is the arrow
$$
\mathcal{ME}(G) \longrightarrow \Map(\Hom_\Z(G,\Q/\Z),\Q/\Z)
$$
that sends $(M,\xi,\psi)$ to the map
$$
\xymatrix{
\Hom_\Z(G,\Q/\Z) \ar[r]^-{(\psi^{-1})^*}_-{\simeq} & \Hom_\Z(H_1(M),\Q/\Z)  
\ar[r]_-\simeq & H_2(M;\Q/\Z) \ar[rr]^-{\phi_{M,\xi}} & &
}
\Q/\Z
$$
where the central map is given by the intersection pairing 
and the last one is the linking quadratic function. 
\end{remark}

\begin{remark}
Originally, Goussarov and Habiro defined finite-type invariants with respect 
to surgery along some kinds of embedded decorated graphs, called \emph{graph claspers}.
By Johnson's result on the generation of the Torelli group \cite{Johnson_generation},
a Torelli surgery can be realized by a finite number of surgeries along graph claspers, and vice-versa.
This implies that the definition of finite-type invariant given here agrees with Goussarov and Habiro's notion.
For integral homology $3$-spheres, the use of the Torelli group and its filtrations to define finite-type invariants
(and, in particular, to recover Ohtsuki's definition \cite{Ohtsuki}) appeared firstly 
in Garoufalidis and Levine's work \cite{GL}.
\end{remark}

\begin{remark}
Some geometric techniques, which are  refered to as \emph{calculus of claspers},
have been developed in \cite{Habiro,GGP} to prove general results about finite-type invariants. 
Calculus of claspers extends\footnote{Calculus of claspers works well when parallelizations (or, equivalently, 
Spin-structures) and homological parametrizations are taken into account \cite[\S 2.2]{Massuyeau_bis}.
We conclude by Lemma \ref{lem:Omega}.} from $\mathcal{M}$ to $\mathcal{ME}(G)$,
so that most of those general results do too.
For instance, it is known \cite{Habiro,Habegger,Garoufalidis_brane} that, for each degree $d$, 
there is only a finite number of linearly independent finite-type invariants for $\mathcal{M}$. 
The statement and the proof given in \cite{Garoufalidis_brane}, or alternatively in \cite{Habegger}, 
apply mutatis mutandis to $\mathcal{ME}(G)$.
\end{remark}

\section{Proof of the main results}
\label{sec:proofs}

In this section, we prove our main results, including Theorem \ref{th:FTI}
and Theorem \ref{th:twist} from the introduction.

\subsection{Computing the RT torsion from a Heegaard splitting}

\label{subsec:computation}

Let $M$ be a closed connected oriented $3$-manifold, 
presented by means of a Heegaard splitting:
$$
M = A \cup B, \quad A\cap B = \partial A = - \partial B,
$$
where $A$ and $B$ are genus $g$ handlebodies.
We recall from \cite{Turaev_spinc} a formula giving the RT
torsion of $M$ equipped with the combinatorial Euler structure that the given Heegaard splitting ``prefers''.
We wish to explore that formula far enough to be able, in the sequel, 
to compare the torsion of $M$ with that of 
the manifold $M_h$ resulting from a Torelli surgery $M\leadsto M_h$.
In particular, we precise which is the geometric Euler structure corresponding 
to that combinatorial Euler structure.\\

\begin{figure}[h]
\centerline{\relabelbox \small 
\epsfxsize 3truein \epsfbox{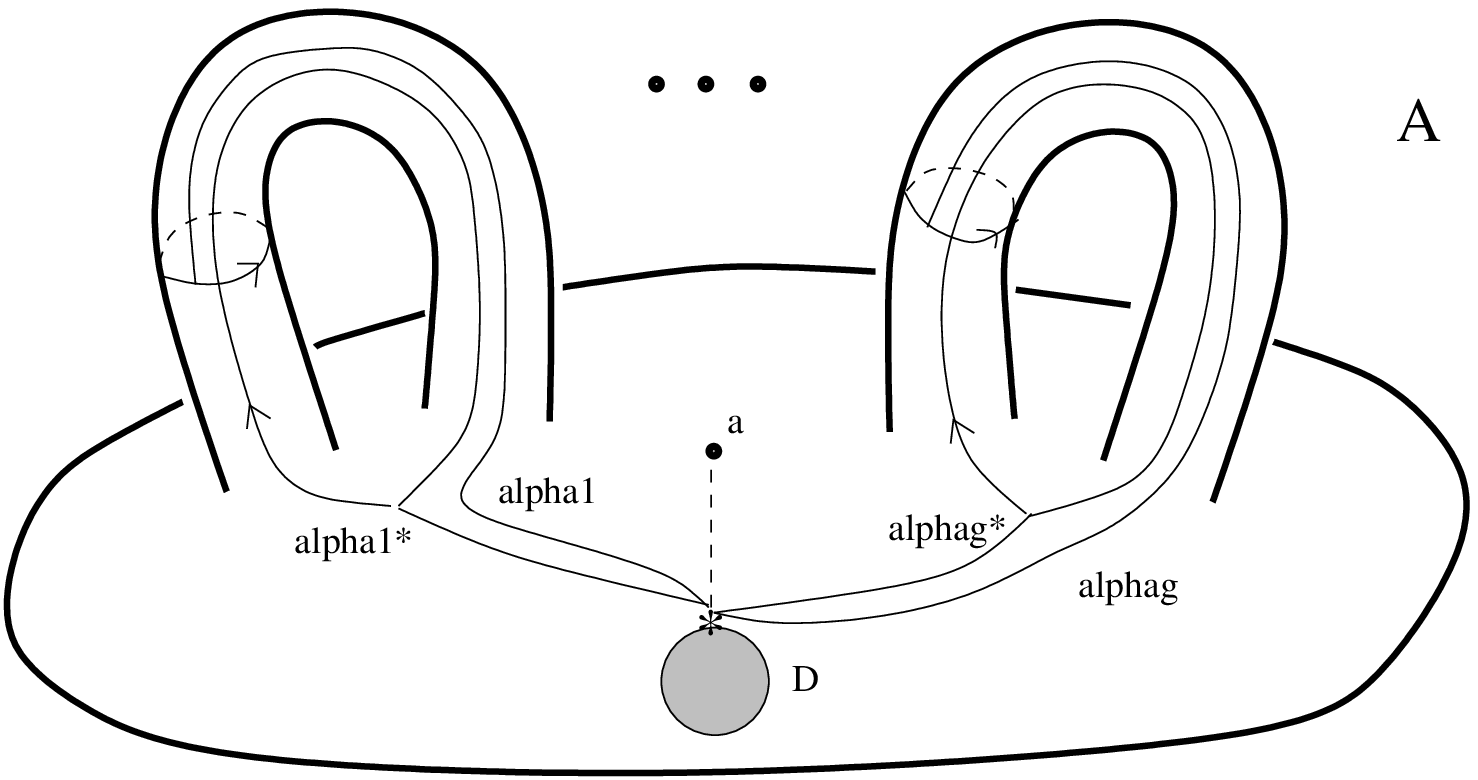}
\adjustrelabel <0.cm,-0.2cm> {D}{$d$}
\adjustrelabel <0.cm,0cm> {A}{$B$}
\adjustrelabel <0.cm,0cm> {a}{$b$}
\adjustrelabel <0cm,0cm> {alpha1}{$\beta_1$}
\adjustrelabel <0.1cm,-0.1cm> {alpha1*}{$\beta_1^*$}
\adjustrelabel <0cm,0cm> {alphag}{$\beta_g$}
\adjustrelabel <0cm,0cm> {alphag*}{$\beta_g^*$}
\adjustrelabel <0.cm,0.1cm> {*}{\large{$\star$}}
\endrelabelbox}

\vspace{0.8cm}

\centerline{\relabelbox \small 
\epsfxsize 3truein \epsfbox{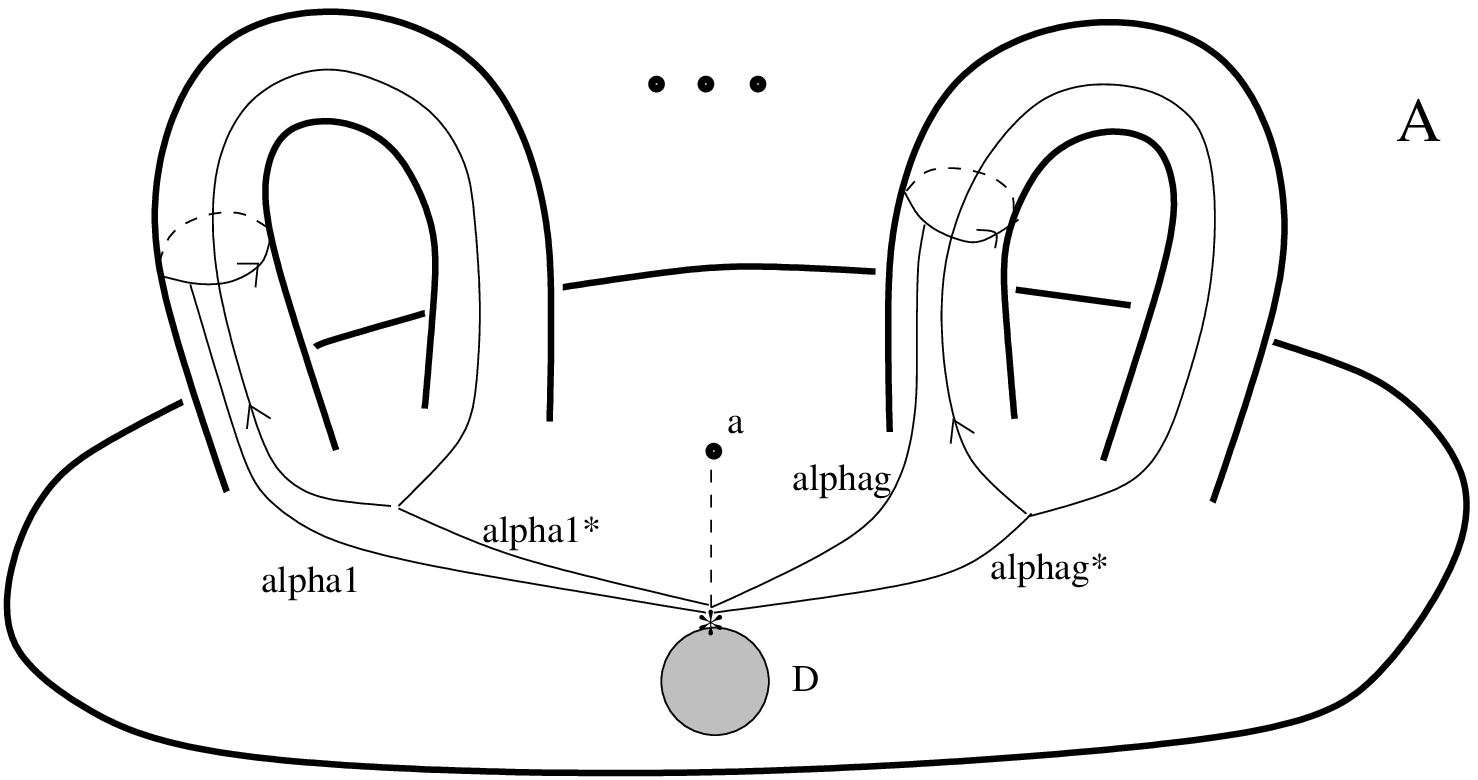}
\adjustrelabel <0.cm,-0.2cm> {D}{$d$}
\adjustrelabel <0.cm,0cm> {A}{$A$}
\adjustrelabel <0.cm,0cm> {a}{$a$}
\adjustrelabel <0.1cm,-0.1cm> {alpha1}{$\alpha_1$}
\adjustrelabel <0.1cm,0.1cm> {alpha1*}{$\alpha_1^*$}
\adjustrelabel <0.1cm,0cm> {alphag}{$\alpha_g$}
\adjustrelabel <0.2cm,0cm> {alphag*}{$\alpha_g^*$}
\adjustrelabel <0.cm,0.1cm> {*}{\large{$\star$}}
\endrelabelbox}
\caption{The lower handlebody $A$ and the upper handlebody $B$ in the Heegaard splitting of $M$.}
\label{fig:Heegaard}
\end{figure}

First of all, we fix some basis of various kinds as indicated on Figure \ref{fig:Heegaard}.
Apart from a common base point $\star$ on the boundary of a common small disk $d\subset A\cap B$,
those choices are either relative to the handlebody $A$, either relative to the handlebody $B$.
Watch out that this picture does not suggest any kind of identification
between $\partial A$ and $\partial B$ elsewhere than in $d$.
In the lower handlebody $A$, the $\alpha_i$'s and the $\alpha_i^*$'s
are respectively meridional and longitudinal loops on $\partial A$. These loops are based at $\star$, 
the basing arc of $\alpha_i$ following the orientation of $\alpha_i^*$.
The point $\star$ itself is connected by a small arc to a point $a$ interior to $A$. 
Similar choices and notations have been fixed for the upper handlebody $B$, 
but the basing arc of $\beta_j$ is now against the orientation  of $\beta_j^*$.

Those choices being fixed, let us consider a cell decomposition $X$
of $M$ induced by the Heegaard splitting:
There is only one $0$-cell $e^0=a$, the center of the ball to which
handles have been added to form $A$; the $1$-cells are $e^1_1,\dots,e^1_g$  
where $e^1_i$ is obtained from the core of the $i$-th handle of $A$,
which is bounded by two points, by adding the trace of those two points 
when the previous ball is ``squeezed'' to $e^0$; 
the $2$-cells are $e^2_1,\dots,e^2_g$ where $e^2_j$ is obtained from the co-core
of the $j$-th handle of $B$, which is bounded by the circle $\beta_j \setminus (\hbox{basing arc})$,
by adding the trace of  that circle in $A$ when $A$ is ``squeezed'' to
$e^0 \cup e^1_1\cup \cdots \cup e^1_g$; there is only one $3$-cell $e^3$, 
namely the complement in $M$ of the cells of smaller dimension. 
As for the orientations, $e^0$ is given the $+$ sign, $e^1_i$ is oriented coherently with
$\alpha_i^*$, $e^2_j$ is oriented so that $e^2_j\bullet \beta_j^*=+1$ and $e^3$ has the orientation of $M$.

Let $p:\widehat{M} \to M$ be the maximal Abelian covering determined by the base point $\star$
and let $\widehat{\star}$ be the distinguished lift of $\star$. The lift of the arc $\xymatrix{\star\ar@{-}[r] & a}$
starting at $\widehat{\star}$ determines a preferred lift $\widehat{e}^0$ of $e_0=a$. 
Similarly, the lift of $\xymatrix{\star\ar@{-}[r] & b}$ starting at $\widehat{\star}$ determines 
a preferred lift $\widehat{b}$ of $b$: Let $\widehat{e}^3$ be the unique lift of $e^3$ containing $\widehat{b}$. 
Let $\widehat{e}^1_i$ be the lift of $e^1_i$ starting at $\widehat{e}^0$. Let $\widehat{e}^2_j$ be the lift of $e^2_j$
contained in $\partial \widehat{e}^3$ with the opposite orientation. Let $\xi$ be the Euler structure
represented by the fundamental family of cells $\widehat{e}$. 
We call it the Euler structure \emph{preferred} by the Heegaard splitting.

The RT torsion of $(M,\xi)$ can be computed
using the cell decomposition $X$. Such a computation involves the sign
$$
\tau_0 := \sgn\ \tau\left(C_*(M \hbox{ cellularized by } X;\R); oo, w\right)
$$
where $oo$ refers to the above choices of $o$rder and $o$rientation of cells,
and $w$ is a basis of $H_*(M;\R)$ representing $\omega_M$.
The boundary operators $\partial_0$ and $\partial_2$ are given by 
$$
\forall i=1,\dots,g, \ \partial \widehat{e}^1_i = (\alpha_i^*-1)\cdot \widehat{e}^0  \quad \textrm{and} \quad
\partial \widehat{e}^3= \sum_{j=1}^g (\beta_j^*-1)\cdot \widehat{e}^2_j.
$$
Thus, the main indeterminate is the boundary operator 
$\partial_1: C_2(\widehat{X}) \to C_1(\widehat{X})$
whose $(i,j)$-minor in the basis $\widehat{e}_{oo}$ is denoted by $\Delta_{ij}$.
The result is as follows. 

\begin{lemma}[Turaev \cite{Turaev_spinc}] 
\label{lem:computation}
Let $\varphi: \Z[H_1(M)] \to \F$ be a ring homomorphism with values in a commutative field 
such that $\varphi(H_1(M))\neq 1$. For any indices $i,j=1,\dots,g$, the following identity holds:
\begin{equation}
\label{eq:Heegaard}
\tau^\varphi(M,\xi)\cdot \left(\varphi(\alpha_i^*)-1\right) \cdot \left(\varphi(\beta_j^*)-1\right) 
= (-1)^{g+i+j+1}\cdot \tau_0 \cdot \varphi(\Delta_{ij}) \in \F.
\end{equation}
\end{lemma}
\noindent
See \cite[Proof of Th. 4.1]{Turaev_spinc} 
or \cite[Proof of Th. II.1.2]{Turaev_bigbook} for the details of  computation.\\

Next, we wish to remove from formula (\ref{eq:Heegaard}) any reference
to the maximal Abelian covering. In particular, we can
find a geometric representant for the Euler structure $\xi$ in the following way.

\begin{figure}[h!]
\centerline{\relabelbox \small 
\epsfxsize 5truein \epsfbox{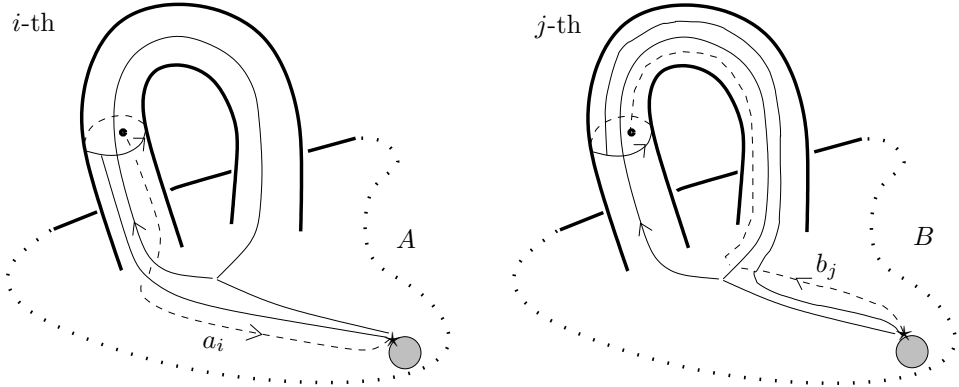}
\adjustrelabel <-0.3cm,0.3cm> {ith}{$i$-th}
\adjustrelabel <-0.3cm,0.3cm> {jth}{$j$-th}
\adjustrelabel <0.cm,0cm> {ai}{$a_i$}
\adjustrelabel <-0.1cm,0.1cm> {bj}{$b_j$}
\adjustrelabel <0.cm,0cm> {A}{$A$}
\adjustrelabel <0.cm,0cm> {B}{$B$}
\adjustrelabel <-0.15cm,-0.15cm> {astar}{\large{$\star$}}
\adjustrelabel <-0.05cm,-0.05cm> {bstar}{\large{$\star$}}
\endrelabelbox}
\caption{The arcs $a_i$'s and $b_j$'s.}
\label{fig:arcs}
\end{figure}

Let $a_i$ be an oriented arc in $A$ connecting 
the center of the $i$-th handle to $\star$ as depicted in Figure \ref{fig:arcs}. 
An oriented arc $b_j$ in $B$, connecting
$\star$ to the center of the $j$-th handle, is depicted on the same figure.
Consider the following Euler chain relative to the cell decomposition $X$:
\begin{equation}
\label{eq:c}
c:= (\xymatrix{\star\ar[r] & a}) + a_1 + \cdots + a_g 
+ b_1 + \cdots + b_g + (\xymatrix{\star & \ar[l] b}).
\end{equation}
It is the image under the covering map $p$ of a ``spider chain'' $\widehat{c}$ 
with body $\widehat{\star}$ and $2g+2$ legs, each connecting
with the appropriate orientation $\widehat{\star}$ to the center of a cell in $\widehat{e}$.
In other words, $c\in \Eul(X)$ is representative for $\xi\in \Eul(M)$. 

Let $f: M \to \R$ be a Morse function, together with a Riemannian metric on $M$, 
such that the Smale condition is satisfied and  the handle decomposition of $M$ induced by $f$ 
is our given Heegaard splitting. More precisely, we assume that $f$ is self-indexing, 
$A=f^{-1}([0,3/2])$ and $B=f^{-1}([3/2,3])$. The function $f$ has $a$ (respectively $b$) 
as only critical point of index $0$ (respectively $3$) and
its critical points of index $1$ (respectively $2$) 
are the centers of the handles of $A$ (respectively $B$).

\begin{figure}[h]
\centerline{\relabelbox \small 
\epsfxsize 3truein \epsfbox{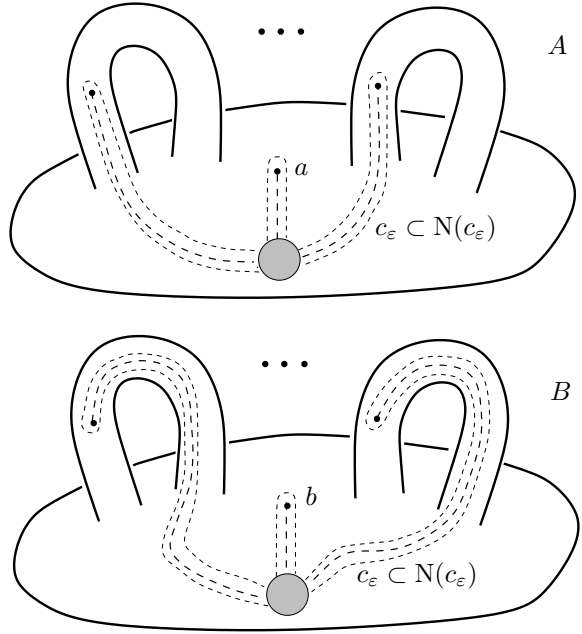}
\adjustrelabel <0.cm,0cm> {A}{$A$}
\adjustrelabel <0.cm,0cm> {B}{$B$}
\adjustrelabel <0.cm,0cm> {a}{$a$}
\adjustrelabel <0.cm,0cm> {b}{$b$}
\adjustrelabel <0.1cm,-0cm> {ca}{$c_\varepsilon \subset \hbox{N}(c_\varepsilon)$}
\adjustrelabel <0cm,-0.1cm> {cb}{$c_\varepsilon \subset \hbox{N}(c_\varepsilon)$}
\endrelabelbox}
\caption{The Euler chain $c_\varepsilon$ representing $\xi$.}
\label{fig:chain}
\end{figure}

\begin{lemma}
\label{lem:Heegaard_smooth}
Let $c_\varepsilon$ be obtained from the chain $c$, 
defined at (\ref{eq:c}), pushing the point $\star$ into the interior of the disk $d$.
Let N$(c_\varepsilon)$ be a ball neighborhood of $c_\varepsilon$
which meets the Heegaard surface $A\cap B$ in the interior of the disk $d$ --
see Figure \ref{fig:chain}.
Then, any non-singular vector field $v$ on $M$ which coincides
with $\nabla f$ outside N$(c_\varepsilon)$, represents the Euler structure $\xi$.
\end{lemma}

\begin{proof}
We have seen that $\xi$ is represented by $c$ as an Euler chain relative to the cell decomposition $X$. 
But, $c$ can also be interpreted as an Euler chain relative to the vector field $\nabla f$. 
We can assume that the cell decomposition $X$
(which has come from the Heegaard splitting)
coincides with the cell decomposition $(X_f,\rho_f)$ that is induced by the Morse function $f$.
By Lemma \ref{lem:Morse}, the combinatorial 
Euler structure $\xi$ corresponds to the geometric Euler structure $\alpha_{\nabla f}([c_\varepsilon])$.
We conclude thanks to the last statement of Lemma \ref{lem:alpha}.
\end{proof}

Finally, the minor $\Delta_{ij}$ in formula (\ref{eq:Heegaard})
can be computed, with no reference to the maximal Abelian covering, thanks to the next lemma.

\begin{lemma} 
\label{lem:minor}
The matrix of the boundary operator $\partial_1: C_2(\widehat{X}) \to C_1(\widehat{X})$ 
in the basis $\widehat{e}_{oo}$ is the image of
$$
\left(\frac{\partial \beta_j}{ \partial \alpha_i^*}\right)_{i,j=1,\dots,g}
$$
under the homomorphism $\incl_*:\pi_1(\partial A\setminus \int\ d, \star) \to H_1(M)$.
Here,  the free derivative $\partial/ \partial \alpha_i^*$ is computed
in the free group $\pi_1(\partial A\setminus \int\ d, \star)$ 
with respect to the basis $(\alpha_i,\alpha_i^*)_{i}$.
\end{lemma}

\begin{proof}
The fact that the boundary operator $\partial_1$ of $\widehat{X}$ can be computed by means
of free differential calculus is classical, but we have to be careful with the basing of the loops.

First, let us recall how free derivatives come into the picture.
Set $\Sigma:=\partial A\setminus \int\ d$ and let $q:\widetilde{\Sigma} \to \Sigma$ be its universal
covering determined by $\star$. Set $\widehat{\Sigma}:= p^{-1}(\Sigma)\subset \widehat{M}$.
For any arc $\gamma \subset \Sigma$ starting at $\star$, 
we denote by $\tilde \gamma$ and $\widehat \gamma$ the lifts of $\gamma$ to $\widetilde{\Sigma}$
and $\widehat{\Sigma}$ respectively, starting at $\tilde \star$ and $\widehat \star$ respectively.
Since $\Sigma$ retracts to a wedge of circles, 
for any loop $\gamma \subset \Sigma$ based at $\star$, the identity
$$
\tilde{\gamma} = \sum_{i=1}^g \frac{\partial  \gamma}{\partial \alpha_i}\cdot \widetilde{\alpha_i}
+\frac{\partial  \gamma}{\partial \alpha_i^*}\cdot \widetilde{\alpha_i^*}
$$
holds in the $\Z[\pi_1\left(\Sigma,\star\right)]$-free module $H_1(\widetilde{\Sigma},q^{-1}(\star))$.
Applying the covering transformation $\widetilde{\Sigma} \to \widehat{\Sigma}$ that sends
$\tilde{\star}$ to $\widehat{\star}$, one gets the identity
\begin{equation}
\label{eq:lift}
\widehat{\gamma} = \sum_{i=1}^g \incl_*\left(\frac{\partial  \gamma}{\partial \alpha_i}\right) \cdot \widehat{\alpha_i}
+\incl_*\left(\frac{\partial  \gamma}{\partial \alpha_i^*}\right)\cdot \widehat{\alpha_i^*}
\end{equation}
in the $\Z[H_1(M)]$-free module $H_1(\widehat{\Sigma},p^{-1}(\star))$. 

We are asked to compute $\partial \widehat{e}^2_j$. 
Let $b'_j$ be the basing arc of $\beta_j$. 
Recall that the $2$-cell $e^2_j$ of $M$ is the core of the $j$-th handle of 
$B$, which we denote by $B_j$, together with the trace of 
$\beta_j\setminus b'_j$ in $A$ during the retraction to its spine. 
(Both $B_j$ and $b'_j$ can be seen on the right hand-side of Figure \ref{fig:arcs}.)
Then, $\beta_j$ is the concatenation of paths 
$b_j' \cdot \partial B_j \cdot b_j'^{ -1}$ (if the loop $\partial B_j$ is based at the endpoint of $b'_j$).
Let $\widehat{B}_j$ be the lift of $B_j$ that is contained in $\widehat{e}^2_j$. 
Then, $\widehat{\beta_j}$ coincides with the loop
$\widehat{b_j'} \cdot \partial \widehat{B}_j \cdot \widehat{b_j'}^{-1}$ (if $\partial \widehat{B}_j$ is based
at the endpoint of $\widehat{b_j'}$). 
The manifold $p^{-1}(A) \subset \widehat{M}$ is a handlebody (in the general sense) 
formed by the lifts of the only $0$-handle
and the $g$ $1$-handles of $A$; it has the $1$-skeleton of $\widehat{X}$ as spine. 
Thus, one sees  that  $\partial \widehat{e}^2_j$ can be computed (as an element of
$C_1(\widehat{X})$) from its intersection with $\partial p^{-1}(A)=\widehat{\Sigma}$. 
More precisely, the cellular $1$-chain $\partial \widehat{e}^2_j$ is obtained from 
$$
\left[\partial \widehat{e}^2_j \cap \widehat{\Sigma}\right] = \left[\partial \widehat{B}_j\right]
= \left[\widehat{\beta_j}\right] \in H_1\left(\widehat{\Sigma},p^{-1}(\star)\right)
$$
by the rules $\widehat{\alpha_i}\mapsto 0$ and $\widehat{\alpha_i^*} \mapsto \widehat{e}^1_i$.
We conclude by applying (\ref{eq:lift}) to $\gamma=\beta_j$.
\end{proof}

\subsection{Finiteness properties of the RT torsion}

\label{subsec:finiteness_properties}

For $G$ a finitely generated Abelian group, we define
$$
\tau: \mathcal{ME}(G) \longrightarrow Q(\Z[G]), \quad 
(M,\xi,\psi) \mapsto Q(\psi^{-1})( \tau(M,\xi)).
$$
In this subsection, some reductions of $\tau$ are shown to be of  finite type
in the sense of Definition \ref{def:fti_structures}.
The study starts with a general result and, next,
splits into two cases: b$_1(M) >0$ and b$_1(M)=0$.

\subsubsection{A general lemma}

Let $G$ be a finitely generated Abelian group.
Instead of the map $\tau: \mathcal{ME}(G) \to Q(\Z[G])$, one can consider the map 
$$
\overline{\tau}: \mathcal{ME}(G) \longrightarrow \Map\left(G\times G,Q(\Z[G])\right)
$$
that sends a triple $(M,\xi,\psi)$ to the pairing $(x,y) \mapsto (x-1) \cdot (y-1) \cdot \tau(M,\xi,\psi)$.
This is the form the torsion has been computed from a Heegaard splitting in \S \ref{subsec:computation}.

\begin{lemma}
\label{lem:equivalence_tau}
The invariants $\overline{\tau}$ and $\tau$ are equivalent.
\end{lemma}
\begin{proof}
Consider two triples $(M_1,\xi_1,\psi_1)$ and $(M_2,\xi_2,\psi_2)$. 
Assuming that $\overline{\tau}(M_1,\xi_1,\psi_1)=\overline{\tau}(M_2,\xi_2,\psi_2)$, that is,
\begin{equation}
\label{eq:zero}
\forall x,y\in G, \quad
\left( \tau(M_1,\xi_1,\psi_1) - \tau(M_2,\xi_2,\psi_2)\right) \cdot (x-1) \cdot (y-1)=0,
\end{equation}
we must prove that $\tau(M_1,\xi_1,\psi_1)=\tau(M_2,\xi_2,\psi_2)$.
If $\rk(G)>0$, we take $x=y\in G$ of infinite order, 
so that $(x-1)$ is invertible in $Q(\Z[G])$.
When $G$ is finite, the difference $\tau(M_1,\xi_1,\psi_1) - \tau(M_2,\xi_2,\psi_2)$, written as
$$
\sum_{g \in G} z(g) \cdot g \in Q(\Z[G]) = \Q [G],
$$
defines a function $z: G \to \Q$. Property (\ref{eq:zero}) means that
$z$ is an affine function $G \to \Q$. So, $z$ must be constant. 
Since $\tau(M_1,\xi_1,\psi_1)$ and $\tau(M_2,\xi_2,\psi_2)$ augment to zero, 
we must have $\sum_{g \in G} z(g)= 0$ and we conclude.
\end{proof}

Let $I$ be the augmentation ideal of the group ring $\Z\left[G\right]$.
Next lemma tells that, for any integer  $d\geq 1$, the reduction of $\overline{\tau}$ 
modulo $\Map\left(G\times G, I^d\right)$ is a finite-type invariant of degree $d-1$.

\begin{lemma}
\label{lem:general}
Let $x,y\in G$. Let $M$ be a connected closed oriented $3$-manifold
with Euler structure $\xi$ and homological parametrization $\psi: G \to H_1(M)$.
For any family $\Gamma$ of $d\geq 1$ Torelli automorphisms 
of the boundaries of pairwise disjoint handlebodies in $M$, we have that
$$
(x-1) \cdot (y-1) \cdot \sum_{\Gamma' \subset \Gamma} 
(-1)^{|\Gamma'|} \cdot 
\tau\left((M,\xi,\psi)_{\Gamma'}\right) \in I^d \
\subset  Q\left(\Z\left[G\right]\right).
$$
\end{lemma}

\begin{proof}

Let $H_1,\dots, H_{d}$ be the handlebodies of the family $\Gamma$,
with corresponding genus $g_1,\dots,g_{d}$ and 
corresponding Torelli automorphisms $h_1,\dots,h_{d}$.
Assume, after isotopy, that each $h_i$ is the identity on a small disk $d_i$. 
Pick a ball in $M$ disjoint from the $H_i$'s,  a disk $D$ on the boundary of this ball,
and connect each $H_i$ to the ball by a $1$-handle attached along $d_i$, at one end,
and outside $D$, at the other end.
Thus, we get a big handlebody $H \subset M$ of genus $g_1 + \cdots + g_{d}$
and we define $h_i: \partial H \to \partial H$ to be the extension of $h_i|_{\partial H_i \setminus d_i}$ by the identity.
For each binary number $\sigma \in \{0,1\}^{d}$, we denote by 
$h^\sigma:\partial H \to \partial H$ the (commuting) product of the $h_i$'s for which $\sigma(i)=1$.
Thus, we are asked to show that
\begin{equation}
\label{eq:goal} 
(x-1)\cdot (y-1)\cdot \sum_{\sigma \in \{0,1\}^{d}} (-1)^{|\sigma|} \cdot  
\tau\left((M,\xi,\psi)_{h^\sigma}\right) \in I^d \ \subset  Q\left(\Z\left[G\right]\right)
\end{equation}
where $|\sigma|$ denotes $\sigma_1+ \cdots + \sigma_{d}$.

Let $H'\subset M$ be a handlebody obtained from $H$ by adding along
$D$ a handle whose core together with $D$ realizes the homology class $\psi(y)$.
Let also $H''$ be a solid torus disjoint from $H'$ and whose core realizes the 
homology class $\psi(x)$. 
Regarding $M\setminus \int(H'\cup H'')$ as a cobordism from $-\partial H'$ to $\partial H''$,
we can find a decomposition of it into handles with
no $0$-handle nor $3$-handle. Furthermore, we can rearrange this cobordism 
and assume that the $1$-handles have been attached along $D$. 
So, we have found a Heegaard splitting 
\begin{equation}
\label{eq:splitting_M_bis}
M=A\cup B
\end{equation}
of $M$ such that the upper handlebody $B$ is obtained from the big handlebody $H$ 
by adding handles along $D$, one of which realizing $\psi(y)$, 
while the lower handlebody $A$ has a handle which realizes the homology class $\psi(x)$.

We apply \S \ref{subsec:computation} to the Heegaard splitting (\ref{eq:splitting_M_bis}).
Recall that this computation begins with various choices (Figure \ref{fig:Heegaard}).
We have choosen the common small disk $d$ and the common base point $\star$ inside the disk $D$,
in such a way that $\partial d\cap \partial D=\star$. 
We have numbered the handles of $B$ in such a way that  the first ones are the handles of $H$
(and this numbering respects the ordering of the handlebodies $H_1,\dots,H_d$)
and the last one realizes $\psi(y)$.
Also, we have ordered the handles of $A$ in such a way that the last one realizes $\psi(x)$.
By equivariance, we can assume that 
the given Euler structure $\xi$ coincides with the Euler structure preferred by the Heegaard splitting.
Then, by Lemma \ref{lem:computation} and Lemma \ref{lem:minor}, 
we have that
\begin{equation}
\label{eq:torsion_of_M_bis}
\tau^\varphi(M,\xi,\psi) \cdot \left(\varphi(x)-1\right) \cdot \left(\varphi(y)-1\right)
= (-1)^{g+1}\cdot \tau_0  \cdot \varphi \psi^{-1} (\Delta_{gg})
\end{equation}
where $\varphi: \Z[G]\to \F$ is a ring homomorphism with values in a commutative field $\F$
and such that $\varphi(G)\neq 1$,  
$\tau_0$ is a certain sign and $\Delta_{gg}$ is the determinant of the matrix
\begin{equation}
\label{eq:minor_M_bis}
\left(\frac{\partial \beta_j}{ \partial \alpha_i^*}\right)_{i,j=1,\dots,g-1}
\end{equation}
after the projection $\incl_*: \Z[\pi_1(\partial A \setminus d,\star)] \to \Z[H_1(M)]$ has been applied.

More generally, for any binary number $\sigma\in \{0,1\}^{d}$, 
the Heegaard splitting (\ref{eq:splitting_M_bis}) of $M$ 
induces a Heegaard splitting of $M_{h^{\sigma}}$:
\begin{equation}
\label{eq:splitting_Mh_bis}
\begin{array}{rcl}
M_{h^\sigma} & = & (M \setminus \int\ H) \cup_{h^{\sigma}} H\\
& = & \left(A \cup \overline{B \setminus H}\right) \cup_{h^\sigma} H\\
& = & A \cup_{\Id \cup h^\sigma} \left(\overline{B \setminus H} \cup H\right)\\
& = & A \cup_{ h_{\hbox{\scriptsize e}}^\sigma } B
\end{array}
\end{equation}
where $h_{\hbox{\scriptsize e}}^\sigma: \partial B \to -\partial A$ denotes the extension 
$\partial B \to \partial B$ of $h^\sigma|_{\partial H \setminus D}$ by the identity,
composed with the equality $\partial B= -\partial A$.
We apply \S \ref{subsec:computation} to this Heegaard splitting of $M_{h^\sigma}$, \emph{keeping} 
the basis that have already been fixed in the case when $\sigma=0\cdots 0$ (Figure \ref{fig:Heegaard}).
In particular, we keep the same common small disk $d$; this makes sense since,
$d$ being included in $D$, it is fixed pointwisely by $h_{\hbox{\scriptsize e}}^\sigma$.

We now prove that the structure $\xi_{h^{\sigma}}\in \Eul(M_{h^\sigma})$ (corresponding
to $\xi\in \Eul(M)$ by the canonical bijection $\Omega_{h^\sigma}$) shifted by the Chillingworth class:
$$
\xi'_{h^\sigma}:=\xi_{h^\sigma} - \incl_*\ t\left( (h^\sigma|_{\partial H\setminus D})/2\right)
$$ 
coincides with the Euler structure preferred by the Heegaard splitting (\ref{eq:splitting_Mh_bis}).
According to Lemma \ref{lem:Heegaard_smooth} applied to $M$,
for a Morse function $f:M\to \R$ inducing the Heegaard splitting (\ref{eq:splitting_M_bis}),
there exists a specific ball $N\subset M$ (which meets the Heegaard surface inside the disk $d$)
and such that any vector field $v$ which coincides with $\nabla f$ outside $N$ represents $\xi$:
We fix one such $v$. The function $f_{h^\sigma_{\hbox{\scriptsize e}}} := f|_{A} \cup_{ h^\sigma_{\hbox{\scriptsize e}} } f|_{B}$
is a Morse function $M_{h^\sigma} \to \R$ which induces the Heegaard splitting (\ref{eq:splitting_Mh_bis}).
Let also $N_{h^\sigma_{\hbox{\scriptsize e}} } \subset M_{h^\sigma}$ be the ball 
$(N\cap A) \cup_{h^\sigma_{\hbox{\scriptsize e}}} (N\cap B)$.
By Lemma \ref{lem:Heegaard_smooth} applied to $M_{h^\sigma}$, the vector field
$v_{h^\sigma_{\hbox{\scriptsize e}} } := v|_{A}  \cup_{h^\sigma_{\hbox{\scriptsize e}}} v|_{B}$ represents the Euler structure 
preferred by the Heegaard splitting (\ref{eq:splitting_Mh_bis}) since it coincides
with $\nabla f_{h^\sigma_{\hbox{\scriptsize e}}}$ outside $N_{h^\sigma_{\hbox{\scriptsize e}}}$. 
Next, since $\nabla f$ is inwards normal to $B$, the vector field $v$ is inwards normal to
$H$ on $\partial H \setminus D$. So, by Lemma \ref{lem:Omega_more_concretely}, 
$\xi_{h^\sigma}'$ is represented by $v_{h^\sigma}:= v|_{M\setminus \int\ H} \cup_{h^\sigma} v|_H$. It remains
now to observe that $v_{h^\sigma}=v_{h^\sigma_{\scriptsize{\hbox{e}}}}$.

Therefore, by Lemma \ref{lem:computation}, we obtain that
\begin{equation}
\label{eq:torsion_of_Mh_bis}
\tau^{\varphi}\left(M_{h^\sigma},\xi_{h^\sigma}',\psi_{h^\sigma}\right) 
\cdot \left(\varphi(x)-1\right) \cdot \left(\varphi (y)-1\right) = 
(-1)^{g+1}\cdot \left(\tau_0\right)_{h^\sigma} 
 \cdot \varphi \psi_{h^\sigma}^{-1}\left(\left(\Delta_{gg}\right)_{h^\sigma}\right),
\end{equation}
where $\varphi: \Z[G]\to \F$ is any ring homomorphism with values in a commutative field $\F$
and such that $\varphi(G)\neq 1$, 
$\left(\tau_0\right)_{h^\sigma}$ is a certain sign while $\left(\Delta_{gg}\right)_{h^\sigma}$ 
is the determinant of the matrix
\begin{equation}
\left(\frac{\partial h^\sigma_{\hbox{\scriptsize e}}(\beta_j)}{ \partial \alpha_i^*}\right)_{i,j=1,\dots,g-1}
\end{equation}
after the projection $\incl_*: \Z[\pi_1(\partial A \setminus d,\star)] 
\to \Z[H_1(M_{h^\sigma})]$ has been applied.

In order to treat the two cases simultaneously (b$_1(M) >0$ and b$_1(M)=0$), 
the following notation from Turaev  will be useful: Let
$\kappa:Q(\Z[G]) \to Q(\Z[G])$ be the group homomorphism 
defined by $\kappa:=\Id$ when $\rk(G)>0$ and 
$$
\kappa(x):= x-\aug(x)\cdot |G|^{-1} \sum_{x\in G} x  \quad \quad \forall x\in Q(\Z[G])=\Q[G]
$$
when $\rk(G)=0$. 
Then, a multiple application of (\ref{eq:torsion_of_Mh_bis}) leads to 
\begin{equation}
\label{eq:torsion_of_Mh_ter}
\tau\left(M_{h^\sigma},\xi_{h^\sigma}',\psi_{h^\sigma}\right) 
\cdot \left(x-1\right) \cdot \left(y-1\right) = 
(-1)^{g+1}\cdot \left(\tau_0\right)_{h^\sigma} 
 \cdot \kappa \psi_{h^\sigma}^{-1}\left(\left(\Delta_{gg}\right)_{h^\sigma}\right).
\end{equation}

Let us prove that all the signs $\left(\tau_0\right)_{h^\sigma}$ are identical. 
Recall from \S \ref{subsec:computation} that $\tau_0$ is the sign of
\begin{equation}
\label{eq:real_torsion_M_bis}
\tau\left( C_*(M \hbox{ cellularized by }X;\R) ; oo, w\right) \in \R\setminus \{0\}
\end{equation}
where $X$ is a cellularization given by the Heegaard splitting (\ref{eq:splitting_M_bis}), 
$oo$ is the cellular basis given by the $o$rder and $o$rientation of cells
and $w$ is a homological basis representing $\omega_M$. 
More generally, $\left(\tau_0\right)_{h^\sigma}$ is the sign of 
\begin{equation}
\label{eq:real_torsion_Mh_bis}
 \tau\left( C_*(M_{h^\sigma} \hbox{ cellularized by }X_{h^\sigma};\R) ; 
oo, w_{h^\sigma} \right) \in \R\setminus \{0\}
\end{equation}
where $X_{h^\sigma}$ is a cellularization given by the Heegaard splitting (\ref{eq:splitting_Mh_bis}),
$oo$ is the cellular basis given by the $o$rder and $o$rientation of its cells
and $w_{h^\sigma}$ is a homological basis representing $\omega_{M_{h^\sigma}}$. 
Consider the graded isomorphism
$$
\iota_*: C_*(M \hbox{ cellularized by }X;\R) \to C_*(M_{h^\sigma} \hbox{ cellularized by }X_{h^\sigma};\R)
$$ 
defined by $oo \mapsto oo$, i.e. obtained by identifying the $o$riented cells of $X$ with those of $X_{h^\sigma}$ 
according to the $o$rderings and dimensions. The map $\iota_*$ is a chain map (in particular,
the identity $\partial _1 \iota_2= \iota_1 \partial_1$ follows from the fact that $h^\sigma$ acts trivially in homology).
At the level of homology, we immediately see that $\iota_0([\star])=[\star]$, 
$\iota_3([M])=[M_{h^\sigma}]$ and $\iota_1([\alpha_i^*])=[\alpha_i^*]$. 
From this latter equality, we deduce that $\iota_1([\beta_j^*])=[\beta_j^*]$ as well (again,
using that $h^\sigma$ preserves the homology). 
So, $\iota_1$ is dual to $\iota_2^{-1}$ with respect to the intersection pairings. Consequently,
$\iota_*$ sends the orientation $\omega_M$ to the orientation $\omega_{M_{h^\sigma}}$ 
and the two torsions (\ref{eq:real_torsion_M_bis}) and (\ref{eq:real_torsion_Mh_bis}) have the same sign.

By taking the alternate sum over $\sigma\in \{0,1\}^{d}$ in (\ref{eq:torsion_of_Mh_ter}), we obtain that
\begin{equation}
\label{eq:alternate_sum}
(x-1)\cdot (y-1)\cdot 
\sum_{\sigma } 
(-1)^{|\sigma|}  \tau\left(M_{h^\sigma},\xi'_{h^\sigma},\psi_{h^\sigma}\right)
= (-1)^{g+1} \tau_0 \cdot 
\kappa\left(\sum_{\sigma } (-1)^{|\sigma|} \psi_{h^\sigma}^{-1}\left(\left(\Delta_{gg}\right)_{h^\sigma}\right)\right).
\end{equation}
From now on, we will not mention anymore the identification $\psi_{h^\sigma}$ 
between the groups $G$ and $H_1(M_{h^\sigma})$.

For each $k=1,\dots,d$, choose a base point $\star_k \in \partial d_k$ and connect it to $\star$
by an arc avoiding $\partial H_l\setminus d_l$ for each $l$.
For each integer $j$ such that $j-(g_1+\cdots+g_{k-1})\in [1,g_k]$ 
(i.e., the $j$-th handle of $B$ comes from the handlebody $H_k$), 
let $\gamma_j \in \pi_1(\partial H_k \setminus d_k,\star_k)'$ be such that 
$h_k$ sends $\beta_j$ to $\gamma_j \cdot \beta_j$ (since it preserves the homology).
We still denote by $\gamma_j$ the element of  $\pi_1(\partial A \setminus d,\star)'$
corresponding to $\gamma_j$ by changing the base points.
It follows that, for any $\sigma \in \{0,1\}^d$ and for any $i=1,\dots,g$,
\begin{equation}
\label{eq:commutator}
\incl_*\left( \frac{\partial h_{\hbox{\footnotesize e}}^\sigma(\beta_j)}{\partial \alpha_i^*}\right) - 
\incl_*\left( \frac{\partial \beta_j}{\partial \alpha_i^*}\right) = 
\left\{\begin{array}{ll}
\incl_*\left( \frac{\partial \gamma_j}{\partial \alpha_i^*}\right)
\in I\subset \Z\left[G\right] & \hbox{ if } \sigma_k=1,\\
0 & \hbox{ if } \sigma_k=0.
\end{array}
\right.
\end{equation}
Then, we get by an induction on $d$ that
\begin{equation}
\label{eq:desalternating}
\sum_{\sigma \in \{0,1\}^{d}} (-1)^{|\sigma|}\cdot (\Delta_{gg})_{h^\sigma}
= (-1)^d \cdot \sum_{s} \Delta_{gg}^s \ \in \Z[G]
\end{equation}
where, on the right handside of the identity, the sum is taken over all 
sub-sequences $s$ of $(1,2,\dots, g_1+\cdots +g_d)$ reaching each of the $d$ intervals 
$$
[1,g_1], g_1+ [1,g_2], \dots, (g_1+ \cdots + g_{d-1}) + [1,g_d]
$$ 
at least one time, and where $\Delta_{gg}^s$ is the determinant of the matrix obtained
from (\ref{eq:minor_M_bis}) by replacing each column whose index $j$ appears in $s$ by
$\left(\partial \gamma_j /\partial \alpha_1^*,\dots, \partial \gamma_j /\partial \alpha_{g-1}^*\right)^t$
and, next, by applying the projection $\incl_*: \Z[\pi_1(\partial A \setminus d,\star)] \to \Z[H_1(M)]$ .
Observe that $\Delta_{gg}^s$ belongs to $I^d$. 
Combining (\ref{eq:alternate_sum}) to (\ref{eq:desalternating}), we obtain that
\begin{equation}
\label{eq:unreduced}
(x-1)\cdot (y-1)\cdot 
\sum_{\sigma\in \{0,1\}^d } 
(-1)^{|\sigma|}  \tau\left(M_{h^\sigma},\xi'_{h^\sigma},\psi_{h^\sigma}\right)
= (-1)^{d+g+1} \tau_0 \cdot  \sum_{s} \Delta_{gg}^s \in I^d.
\end{equation}
 
Applying the next statement to 
$z_k:=\incl_*\ t\left(h_k|_{\partial H_k \setminus d_k}\right)/2$ 
for $k=1,\dots,d$, we obtain (\ref{eq:goal}) which we aimed to prove.

\begin{claim} For any $z_1,\dots,z_d \in G$, the alternate sum
$$
(x-1) \cdot (y-1)\cdot 
\sum_{\sigma \in \{0,1\}^{d}} (-1)^{|\sigma|}  \cdot 
\tau(M_{h^\sigma},\xi'_{h^\sigma},\psi_{h^\sigma})
 \cdot\prod_{k|\sigma_k=1} z_k
$$
is equal modulo $I^{d+1}$ to 
$$ 
(-1)^{d+g+1} \tau_0 \cdot \sum_{\sigma \in \{0,1\}^{d}} \mathcal{D}_{\{k|\sigma_k=1\}} \cdot
\prod_{k|\sigma_k=0} (z_k-1) \in I^d/I^{d+1}.
$$
Here, for any $S \subset \{1,\dots,d\}$,
$\mathcal{D}_S$ denotes the sum of the determinants of all the matrices obtained from
$$
\left(\aug\left( \frac{\partial \beta_j}{ \partial \alpha_i^*}\right)\right)_{i,j=1,\dots,g-1}
$$
by replacing, for each $k\in S$, exactly one column of index $j \in (g_1+\cdots+g_{k-1}) + [1,g_k]$ 
by the column vector $\left(\incl_*(\partial \gamma_{j} /\partial \alpha_1^*),\dots, 
\incl_*(\partial \gamma_{j} /\partial \alpha_{g-1}^*)\right)^t$.
\end{claim}

This claim is proved by a double induction on $(d,n)$ where $n:=|\{k|z_k\neq 1\}|\leq d$.
When $n=0$, the claim is obtained by reducing (\ref{eq:unreduced}) modulo $I^{d+1}$ 
since, to relate some above two notations, we have that
$$
\sum_{s} \Delta_{gg}^s  = 
\sum_{k=1}^d \ \sum_{j_k=g_1 + \cdots + g_{k-1}+1}^{g_1 + \cdots + g_{k-1}+g_k}
 \Delta_{gg}^{(j_1,\dots,j_d)} = \mathcal{D}_{\{1,\dots,d\}} \ \in I^d/I^{d+1}.
$$
When $d=1$, we have modulo $I^2$ that
\begin{eqnarray*}
&& (x-1) \cdot (y-1) \cdot \left(\tau\left(M,\xi,\psi\right) -
\tau\left(M_{h_1},\xi_{h_1}',\psi_{h_1}\right) \cdot z_1\right)\\
&=&  (x-1) \cdot (y-1) \cdot \left( \left(\tau\left(M,\xi,\psi\right) -
\tau\left(M_{h_1},\xi'_{h_1},\psi_{h_1}\right)\right) 
- (z_1-1) \cdot \tau\left(M_{h_1},\xi'_{h_1},\psi_{h_1}\right)\right)\\
&=& (-1)^{1+g+1} \tau_0 \sum_{j_1=1}^{g_1} \Delta_{gg}^{(j_1)}
- (z_1-1) \cdot (-1)^{g+1}\tau_0\cdot  \kappa\left((\Delta_{gg})_{h_1}\right)\\
&=& (-1)^{1+g+1} \tau_0 \cdot  \left(\sum_{j_1=1}^{g_1} \Delta_{gg}^{(j_1)} +
(z_1-1)\cdot \aug (\Delta_{gg}) \right)\\
&=& (-1)^{1+g+1} \tau_0 \cdot \left( \mathcal{D}_{\{1\}} \cdot 1 + \mathcal{D}_{\varnothing} \cdot (z_1-1)  \right)
\end{eqnarray*}
where the second equality follows from (\ref{eq:torsion_of_Mh_ter}) and (\ref{eq:unreduced}),
while the third equality follows from the facts that $\kappa(u)\cdot (v-1)=u\cdot (v-1)$, 
for all $u\in \Z[G]$ and $v\in G$, and that $\aug(\Delta_{gg}) = \aug\left((\Delta_{gg})_{h_1}\right)$.

Assume now that $d > 1$ and $n>0$ and that the claim holds at $(d-1,n-1)$ and at $(d,n-1)$. 
Up to re-numbering of the $z_k$'s, we can assume that $z_d\neq 1$. We have
\begin{eqnarray*}
&& (x-1) \cdot (y-1) \cdot \sum_{\sigma \in \{0,1\}^{d}} (-1)^{|\sigma|} \cdot 
\tau\left(M_{h^\sigma},\xi'_{h^\sigma},\psi_{h^\sigma}\right)  \cdot\prod_{k|\sigma_k=1} z_k\\
&=& (x-1)\cdot (y-1) \cdot \Bigg( \left(\sum_{\sigma \in \{0,1\}^{d}} (-1)^{|\sigma|} \cdot
\tau\left(M_{h^\sigma},\xi'_{h^\sigma},\psi_{h^\sigma}\right)  \cdot\prod_{k|\sigma_k=1} z_k\right)_{z_d=1}\\
&& - (z_d-1) \sum_{\sigma \in \{0,1\}^{d-1}} (-1)^{|\sigma|}  
\tau\left((M_{h_d})_{h^\sigma},(\xi'_{h_d})'_{h^\sigma},(\psi_{h_d})_{h^\sigma}\right) \cdot\prod_{k|\sigma_k=1} z_k \Bigg).
\end{eqnarray*}
We apply the induction hypothesis to the above two summands.
The first one is of type $(d,n-1)$. The second one is of type $(d-1,n-1)$ 
and relative to $(M_{h_d},\xi'_{h_d},\psi_{h_d})$,
but note that $\aug( \partial \beta_j/\partial \alpha_i^*) = \aug( \partial h_d(\beta_j)/\partial \alpha_i^*)$.
Finally, we get
\begin{eqnarray*}
&&(x-1) \cdot (y-1) \cdot \sum_{\sigma \in \{0,1\}^{d}} (-1)^{|\sigma|} \cdot 
\tau\left(M_{h^\sigma},\xi'_{h^\sigma},\psi_{h^\sigma}\right)  \cdot\prod_{k|\sigma_k=1} z_k \\
&=& \left((-1)^{d+g+1} \tau_0 \cdot \sum_{\sigma \in \{0,1\}^{d}} \mathcal{D}_{\{k|\sigma_k=1\}} \cdot
\prod_{k|\sigma_k=0} (z_k-1)\right)_{z_d=1}\\
&& -(z_d-1)\cdot (-1)^{(d-1)+g+1} \tau_0 \cdot \sum_{\sigma \in \{0,1\}^{d-1}} \mathcal{D}_{\{k|\sigma_k=1\}} \cdot
\prod_{k|\sigma_k=0} (z_k-1)\\
&=&(-1)^{d+g+1} \tau_0 \cdot \sum_{\sigma \in \{0,1\}^{d}} 
\mathcal{D}_{\{k|\sigma_k=1\}} \cdot \prod_{k|\sigma_k=0} (z_k-1).
\end{eqnarray*}
\end{proof}

\subsubsection{Case of manifolds $M$ with b$_1(M)>0$}

\label{subsubsec:finiteness_properties_nonQHS}

\begin{lemma}[\cite{Turaev_bigbook}]
\label{lem:back}
Let $G$ be a finitely generated Abelian group of positive rank,
and let $p\in G$ be primitive in $G/\Tors\ G$. If $I$ denotes the augmentation ideal of $\Z[G]$, 
then we have that
$$
\forall x\in Q(\Z[G]), \forall k\geq 0, \  
x(p-1) \in I^{k+1}\Longrightarrow x\in I^k.
$$
\end{lemma}

\begin{proof}
The proof is similar to that of \cite[Lemma II.2.5]{Turaev_bigbook}.

\begin{claim}
Let $G'$ be a direct summand of $\langle p \rangle$ in $G$. Then, we have
$$
I^{k+1} = (p-1) \cdot I^k + I^{k+1} \cap \Z[G'].
$$
\end{claim}

Assuming this, we see that there exist $a\in I^k$ and $b \in I^{k+1} \cap \Z[G']$ such that
$x(p-1)= a(p-1) + b$. Projecting this identity in $\Z[G']$, we get $b=0$ and $x=a$
since, $p$ being of infinite order, $p-1$ is invertible in $Q(\Z[G])$.

Let $g_1,\dots, g_{r-1}$ be some generators of $G'$ and set $g_r:=p$. The ideal $I$ is additively generated
by its elements of the type $(g_i-1)y$, where $y\in G$ and $i=1,\dots,r$. So, $I^{k+1}$ is
additively generated by the products $y_1\cdots y_{k+1} \cdot y$, 
where each $y_j$ is of the form  $(g_i-1)$ and $y\in G$. If such a product
$z:=y_1\cdots y_{k+1} \cdot y$ has no factor $(p-1)$ and if $y$ writes as
$p^n \cdot y'$ where $y'\in G'$ and $n\neq 0$,  then one can apply one of the two identities
$p=(p-1)+1$ and $p^{-1}=-p^{-1}\cdot (p-1)+1$ in order to reduce $|n|$. 
So, we conclude that $z$ belongs to  $(p-1) \cdot I^k + I^{k+1} \cap \Z[G']$.
\end{proof}

\begin{theorem}
\label{th:nonQHS}
Let $G$ be a finitely generated Abelian group of positive rank.
Let $M$ be a connected closed oriented $3$-manifold
with Euler structure $\xi$ and homological parametrization $\psi: G \to H_1(M)$.
For any family $\Gamma$ of $d\geq 2$ Torelli automorphisms 
of the boundaries of pairwise disjoint handlebodies in $M$, we have that
$$
 \sum_{\Gamma' \subset \Gamma} 
(-1)^{|\Gamma'|} \cdot 
\tau\left((M,\xi,\psi)_{\Gamma'}\right) \in I^{d-2} \
\subset  Q\left(\Z\left[G\right]\right)
$$
where $I$ is the augmentation ideal of the group ring $\Z\left[G\right]$.
\end{theorem}

\begin{proof} Apply Lemma \ref{lem:general} taking $x=y$ to be primitive
in $G/ \Tors\ G$ and $d\geq 2$. Next, apply Lemma \ref{lem:back} two times.
\end{proof}

Recall that $\tau(M,\xi) \in \Z\left[H_1(M)\right]$ when b$_1(M) >1$. 
When b$_1(M)=1$, this is not true anymore but 
$\tau(M,\xi)$ has a well-defined \emph{polynomial part} \cite[\S II.3]{Turaev_bigbook}.
More precisely, let $t \in H_1(M)$ be a generator up to torsion.
Let $K_t(\xi) \in \Z$ be such that 
$$
c(\xi)- K_t(\xi)\cdot t \in \Tors\ H_1(M), 
$$
where $c(\xi) \in H_1(M)$ denotes the Chern class.\footnote{Recall that,
for $\xi=[v]\in \Eul(M)$, $c(\xi)\in H^2(M) \simeq H_1(M)$ 
is the obstruction to find a non-singular section of T$M/\langle v\rangle$.} Set
$$
[\tau](M,\xi) := \tau(M,\xi) + t^{K_t(\xi)/2} \cdot
\frac{\sum_{x \in \Tors\ H_1(M)} x}{(t-1)\cdot(t^{-1}-1)}
\in Q(\Z[H_1(M)]).
$$
Then, $[\tau](M,\xi)$ does not depend on the choice of $t$
and belongs to $\Z[H_1(M)]$. Since the canonical bijections 
$\Omega_h:\Eul(M) \to \Eul(M_h)$ and $\Phi_h: H_1(M) \to H_1(M_h)$
induced by a Torelli surgery $M\leadsto M_h$ commute with the Chern class maps, 
we deduce the following 

\begin{corollary}
\label{cor:rank_one}
Let $G$ be a finitely generated Abelian group of rank one.
Let $M$ be a connected closed oriented $3$-manifold
with Euler structure $\xi$ and homological parametrization $\psi: G \to H_1(M)$.
For any family $\Gamma$ of $d\geq 2$ Torelli automorphisms 
of the boundaries of pairwise disjoint handlebodies in $M$, we have that
$$
 \sum_{\Gamma' \subset \Gamma} 
(-1)^{|\Gamma'|} \cdot 
[\tau]\left((M,\xi,\psi)_{\Gamma'}\right) \in I^{d-2} \
\subset  \Z\left[G\right]
$$
where $I$ is the augmentation ideal of the group ring $\Z\left[G\right]$.
\end{corollary}

\subsubsection{Case of manifolds $M$ with b$_1(M)=0$}

\label{subsubsec:finiteness_properties_QHS}

In the case of rational homology spheres, there does not seem to exist
immediate analogue of Lemma \ref{lem:back} except in the cyclic case.

\begin{lemma}
\label{lem:back_QHS}
Let $G$ be a finite cyclic group with generator $p$.
Let $\kappa: \Q[G] \to \Q[G]$ be defined by $\kappa(x):= x-\aug(x)\cdot |G|^{-1} \Sigma$
where $\Sigma:=\sum_{x\in G} x$.
If $I$ denotes the augmentation ideal of $\Z[G]$, then we have that
$$
\forall x\in \Q[G], \forall k\geq 0, \  
x(p-1) \in I^{k+1}\Longrightarrow \kappa(x) \in \kappa(I^k).
$$
\end{lemma}

\begin{proof} Since $I^{k+1}=(p-1)\cdot I^k$, we find $a\in I^k$ such that $x(p-1)=a(p-1)$.
Using the identity
$$
(p-1)\cdot \left(1+2p+3p^2+\cdots+ m p^{m-1}\right)= m - \Sigma
$$
where $m:=|G|$, one gets $x(m - \Sigma)= a(m-\Sigma)$, i.e. $\kappa(x)=\kappa(a)$.
\end{proof}

\begin{theorem}
\label{th:QHS}
Let $G$ be a finite cyclic group.
Let $M$ be a connected closed oriented $3$-manifold
with Euler structure $\xi$ and homological parametrization $\psi: G \to H_1(M)$.
For any family $\Gamma$ of $d\geq 2$ Torelli automorphisms 
of the boundaries of pairwise disjoint handlebodies in $M$, we have that
$$
 \sum_{\Gamma' \subset \Gamma} 
(-1)^{|\Gamma'|} \cdot 
\tau\left((M,\xi,\psi)_{\Gamma'}\right) \in \kappa(I^{d-2}) \
\subset  Q(\Z\left[G\right])=\Q[G]
$$
where $I$ is the augmentation ideal of the group ring $\Z\left[G\right]$.
\end{theorem}

\begin{proof}
Apply Lemma \ref{lem:general} taking $x=y=p$ generating 
$G$ and $d\geq 2$. Observe that 
$\kappa(I^k)=I^k$ for $k>0$ and recall that $\aug(\tau(M,\xi))=0$.
Apply Lemma \ref{lem:back_QHS} two times.
\end{proof}

\subsection{Non-variation of the RT torsion}

\label{subsec:non-variation}

We now give sufficient conditions on the instructions of a Torelli surgery
assuring that the RT torsion does not change.

\begin{lemma} 
\label{lem:non-variation}
Let $G$ be a finitely generated Abelian group.
Let $M$ be a closed oriented connected $3$-manifold with Euler structure $\xi$ 
and homology parametrized by $\psi:G \to H_1(M)$.
Let $H$ be a handlebody in $M$ and let $h$ be a Torelli automorphism of $\partial H$
which is the identity on a disk $D \subset \partial H$.
If $h$ acts trivially on $\pi_1(\partial H \setminus D,\star)$ modulo $\eta'=[\eta,\eta]$ where
$\star \in \partial D$ and 
$$
\eta:= \Ker\left( \xymatrix{\pi_1(\partial H \setminus D, \star) \ar@{->>}[r]
&H_1(\partial H \setminus D) \ar[r]^-{\incl_*} & H_1(M) }\right),
$$
then we have that
$$
\psi^{-1} \circ \incl_* \left( t(h|_{\partial H\setminus D})/2 \right) \cdot
\tau(M,\xi,\psi)= \tau\left((M,\xi,\psi)_h\right) \in Q\left(\Z[G]\right)
$$
where $t$ denotes the Chillingworth homomorphism.
\end{lemma}

This lemma has the following immediate two consequences.

\begin{theorem}
\label{th:null-claspers}
Let $G$ be a finitely generated Abelian group.
Let $M$ be a closed oriented connected $3$-manifold with Euler structure $\xi$ 
and homology parametrized by $\psi:G \to H_1(M)$.
Let $H$ be a handlebody in $M$ and let $h$ be a Torelli automorphism of $\partial H$.
If the handles of $H$ are null-homologous in $M$, then we have that
$$
\tau(M,\xi,\psi)= \tau\left((M,\xi,\psi)_h\right) \in Q\left(\Z[G]\right).
$$
\end{theorem}

\begin{theorem}
\label{th:solvable}
Let $G$ be a finitely generated Abelian group.
Let $M$ be a closed oriented connected $3$-manifold with Euler structure $\xi$
and homology parametrized by $\psi: G \to H_1(M)$.
Let $H$ be a handlebody in $M$ and let $h$ be a Torelli automorphism of $\partial H$
which is the identity on a disk $D \subset \partial H$.
If $h$ acts trivially on the second solvable quotient $\pi/\pi''$
of $\pi:= \pi_1(\partial H\setminus D,\star)$ where $\star \in \partial D$, then
$$
\tau(M,\xi,\psi)= \tau\left((M,\xi,\psi)_h\right) \in Q\left(\Z[G]\right).
$$
\end{theorem}

As for the latter theorem, we should comment that Torelli automorphisms
of $\Sigma_{g,1}$ (:= the connected genus $g$ surface with $1$ boundary component)
which act trivially on the second solvable quotient of its fundamental group $\pi$, do exit.
For this, let us recall  the Magnus representation of the Torelli group of $\Sigma_{g,1}$
$$
r^{\mathfrak{\scriptsize a}} :
\mathcal{T}\left(\Sigma_{g,1}\right) \longrightarrow \hbox{GL}\left(2g;\Z[\pi/\pi']\right), \quad
h \longmapsto  \mathfrak{a}
\left(\overline{\frac{\partial h_*(z_j)}{\partial z_i}}\right)_{i,j}
$$
introduced in Remark \ref{rem:Chillingworth}. In general,
for any normal subgroup $\eta$ of $\pi$ and for any word $w\in \pi$, one has that
\begin{equation}
\label{eq:free_calculus}
w \in \eta' \Longleftrightarrow \forall i \in \{1,...,2n\},
\mathfrak{q}\left(\frac{\partial w}{\partial z_i}\right)=0 \in \Z[\pi/\eta]
\end{equation}
where $\mathfrak{q}$ is the ring homomorphism induced by the quotient map $\pi \to \pi/\eta$.
Taking $\eta=\pi'$, one sees that the kernel of $r^{\mathfrak{\scriptsize a}}$
consists of all the Torelli automorphisms which act trivially on the second solvable quotient.
Elements of that kernel have been exhibited by Suzuki \cite{Suzuki, Suzuki_bis}
(see also \cite{Perron}). On the other hand,
one sees using Remark \ref{rem:Chillingworth},
that the Chillingworth homomorphism vanishes on such diffeomorphisms.

\begin{proof}[Proof of Lemma \ref{lem:non-variation}]

It suffices to prove that, for any ring homomorphism $\varphi$ from the group ring
$\Z[G]$ to a commutative field $\F$,
\begin{equation}
\label{eq:reduction}
\tau^\varphi (M,\xi,\psi)= \tau^{\varphi}(M_h,\xi_h - \incl_*\ t(h|_{\partial H\setminus D})/2,\psi_h) \in \F.
\end{equation}
When $\varphi(G)=1$, this identity is trivial.
So, we assume that there exists a $x\in G$ such that $\varphi(x)\neq 1$.

We follow the proof of Lemma \ref{lem:general} from which we keep the notations. 
There are fewer variables since we are working now with a single handlebody $H_1=H$ (i.e, $d=1$)
and we have taken $y=x$. From that proof, we retain the two formulas (\ref{eq:torsion_of_M_bis})
and (\ref{eq:torsion_of_Mh_bis}):
\begin{equation}
\label{eq:torsion_of_M}
\tau^\varphi(M,\xi,\psi) \cdot \left(\varphi(x)-1\right)^2
= (-1)^{g+1}\cdot \tau_0  \cdot \varphi \psi^{-1} (\Delta_{gg})
\end{equation}
\begin{equation}
\label{eq:torsion_of_Mh}
\tau^\varphi(M_h,\xi_h',\psi_h) \cdot \left(\varphi(x)-1\right)^2
= (-1)^{g+1}\cdot (\tau_0)_h  \cdot \varphi  \psi_h^{-1} \left(\left(\Delta_{gg}\right)_h\right).
\end{equation}
Recall that $\xi'_h$ denotes the shift of $\xi_h$ caused by the Chillingworth class, that
the two signs $\tau_0$ and $(\tau_0)_h$ are identical, and that
$$
\Delta_{gg} = \det
\left(\incl_*\left( \frac{\partial \beta_j}{ \partial \alpha_i^*} \right)\right)_{i,j=1,\dots,g-1}
\quad \hbox{vs} \quad
\left(\Delta_{gg}\right)_h = \det
\left(\incl_*\left( \frac{\partial h_{\hbox{\scriptsize e}}(\beta_j)}{ \partial \alpha_i^*} \right) \right)_{i,j=1,\dots,g-1}.
$$

We now use the extra assumption on $h$.
For any $x\in\pi_1(\partial H \setminus D,\star)$, the class 
$h_*(x)$ differs from $x$ by something in $\eta'$.
Since $h_{\hbox{\scriptsize e}}$ is the extension of $h$ by the identity,
$h_{\hbox{\scriptsize e}}(\beta_j)$ differs from $\beta_j$ 
in $\pi_{\hbox{\scriptsize e}}:= \pi_1(\partial A \setminus d,\star)$ by something in ${\eta_{\hbox{\scriptsize e}}}'$, where
$$
\eta_{\hbox{\scriptsize e}}:= \Ker\left( \xymatrix{\pi_1(\partial A \setminus d,\star)  \ar[r]^-{\incl_*} & H_1(M) }\right).
$$
We deduce from (\ref{eq:free_calculus}) that
$$
\incl_*\left(\frac{\partial h_{\hbox{\scriptsize e}}(\beta_j)}{ \partial \alpha_i^*}\right)
= \incl_*\left(\frac{\partial \beta_j}{ \partial \alpha_i^*}\right)
\in \Z\left[\pi_{\hbox{\scriptsize e}}/\eta_{\hbox{\scriptsize e}}\right]= \Z\left[H_1(M)\right].
$$
We conclude that the above two determinants are essentially equal and, so, that
$$
\tau^\varphi(M,\xi,\psi) = \tau^\varphi(M_h,\xi'_h,\psi_h) \in \F. 
$$
\end{proof}

\subsection{Going deeper into the Johnson filtration}

\label{subsec:Johnson}

The results from \S \ref{subsec:finiteness_properties} can be improved
if one restricts the surgery to some subgroups of the Torelli group.

Let $ \mathcal{M}\left( \Sigma_{g,1} \right)$ denote the
mapping class group of $\Sigma_{g,1}$. For each $c\geq 0$, 
let $\mathcal{M}\left(\Sigma_{g,1}\right)[c]$ be the subgroup 
containing those automorphisms acting trivially on the $c$-th nilpotent quotient $\pi/\pi_{c+1}$ 
of $\pi:=\pi_1\left(\Sigma_{g,1},\star\right)$. This leads to the \emph{Johnson filtration}
$$
\mathcal{M}\left(\Sigma_{g,1}\right) = \mathcal{M}\left(\Sigma_{g,1}\right)[0]
\supset \mathcal{M}\left(\Sigma_{g,1}\right)[1] \supset \mathcal{M}\left(\Sigma_{g,1}\right)[2] \supset \cdots
$$
An automorphism of $\Sigma_g$ (:= the closed oriented connected
surface of genus $g$) is said to be of \emph{class} $c$, if it arises from $\mathcal{M}\left(\Sigma_{g,1}\right)[c]$
through the inclusion $\Sigma_{g,1} \subset \Sigma_g$. For instance, the automorphisms of class $1$ are the Torelli automorphisms. 

\begin{lemma}
\label{lem:general_Johnson}
Let $G$ be a finitely generated Abelian group and let $x,y\in G$.
Let $M$ be a connected closed oriented $3$-manifold
with Euler structure $\xi$ and homological parametrization $\psi: G \to H_1(M)$.
For any family $\Gamma$ of $d\geq 1$ Torelli automorphisms having class $c\geq 1$
of the boundaries of pairwise disjoint handlebodies in $M$, we have that
$$
(x-1) \cdot (y-1) \cdot \sum_{\Gamma' \subset \Gamma} 
(-1)^{|\Gamma'|} \cdot 
\tau\left((M,\xi,\psi)_{\Gamma'}\right) \in I^{c \cdot d} \
\subset  Q\left(\Z\left[G\right]\right)
$$
where $I$ is the augmentation ideal of the group ring $\Z\left[G\right]$.
\end{lemma}

\begin{proof}
The proof is exactly the same as for Lemma \ref{lem:general}
except that, at the identity (\ref{eq:commutator}), we can be more precise:
$$
\incl_*\left( \frac{\partial h_{\hbox{\footnotesize e}}^\sigma(\beta_j)}{\partial \alpha_i^*}\right) - 
\incl_*\left( \frac{\partial \beta_j}{\partial \alpha_i^*}\right) = 
\left\{\begin{array}{ll}
\incl_*\left( \frac{\partial \gamma_j}{\partial \alpha_i^*}\right)
\in I^c \subset \Z\left[G\right] & \hbox{ if } \sigma_k=1,\\
0 & \hbox{ if } \sigma_k=0,
\end{array}\right.
$$
because $\gamma_j$ belongs now to $\pi_1(\partial A \setminus d,\star)_{c+1}$.
Thus, we obtain in place of (\ref{eq:unreduced})
$$
(x-1)\cdot (y-1)\cdot 
\sum_{\sigma\in \{0,1\}^d } 
(-1)^{|\sigma|}  \tau\left(M_{h^\sigma},\xi'_{h^\sigma},\psi_{h^\sigma}\right)
= (-1)^{d+g+1} \tau_0 \cdot  \sum_{s} \Delta_{gg}^s \in I^{cd}.
$$
If $c=1$, we repeat the proof of Lemma \ref{lem:general}. If $c> 1$, the demonstration stops here, 
since the Chillingworth class of each gluing diffeomorphism $h_k$ is then trivial 
\cite[\S 5]{Johnson_homomorphism}.
\end{proof}

Next, the results from \S \ref{subsubsec:finiteness_properties_nonQHS} and
\S \ref{subsubsec:finiteness_properties_QHS} extend in the obvious way to any class $c \geq 2$. 
For instance, the analogue of Theorem \ref{th:nonQHS} is

\begin{theorem}
Let $G$ be a finitely generated Abelian group of positive rank.
Let $M$ be a connected closed oriented $3$-manifold
with Euler structure $\xi$ and homological parametrization $\psi: G \to H_1(M) $.
For any family $\Gamma$ of $d\geq 1$ Torelli automorphisms having class $c\geq 2$
of the boundaries of pairwise disjoint handlebodies in $M$, we have that
$$
 \sum_{\Gamma' \subset \Gamma} 
(-1)^{|\Gamma'|} \cdot 
\tau\left((M,\xi,\psi)_{\Gamma'}\right) \in I^{cd-2} \
\subset  Q\left(\Z\left[G\right]\right)
$$
where $I$ is the augmentation ideal of the group ring $\Z\left[G\right]$.
\end{theorem}

\begin{example}
In particular, taking $d=1$, we see that a Torelli surgery $M\leadsto M_h$ 
whose gluing diffeomorphism $h$ is of class $c\geq 2$,
can not change the RT torsion more than by an element of $I^{c-2}$.
\end{example}

\section{Applications and examples}
\label{sec:applications}

As a conclusion, we give some applications of the finiteness properties for
the Reidem\-eis\-ter--Turaev torsion.

\subsection{Recovering some well-known finiteness properties}

\label{subsec:check}

\subsubsection{The Casson--Walker--Lescop invariant}

Recall from \cite[\S VII.4]{Turaev_bigbook} that the Reide\-meister--Turaev torsion determines
the Casson--Walker--Lescop invariant of a closed oriented connected $3$-manifold $M$ with b$_1(M)>0$:
$$
\lambda(M)= \left\{ 
\begin{array}{ll}
\aug\left([\tau](M,\xi)\right) - |\Tors\ H_1(M)|/12 & \hbox{ if } \hbox{b}_1(M)=1,\\
(-1)^{\hbox{\footnotesize b}_1(M)+1} \aug\left(\tau(M,\xi)\right) & \hbox{ if } \hbox{b}_1(M)>1.
\end{array}\right.
$$
Here, $[\tau](M,\xi)$ denotes the polynomial part of $\tau(M,\xi)$ whose definition 
has been recalled at \S \ref{subsubsec:finiteness_properties_nonQHS}.

Applying Theorem \ref{th:nonQHS}  and Corollary \ref{cor:rank_one} with $d=3$, 
we recover the well-known fact that the Casson--Walker--Lescop invariant is a degree $\leq 2$ invariant.
This also follows from Lescop's sum formula \cite{Lescop}. 

\subsubsection{The Alexander--Conway polynomial}

Let $M$ be a closed oriented connected $3$-manifold with b$_1(M)=1$.
Recall that the Alexander--Conway polynomial of $M$ is defined by
$$
\nabla_M(z) :=|\Tors\ H_1(M)|^{-1} \cdot \left. \Delta_M(t) \right|_{z^2=t+t^{-1}-2} \in \Q\left[z^2\right],
$$
where $\Delta_M(t) \in \Z[t^{\pm 1}]$ denotes the Alexander polynomial of $M$ 
symmetrized so that $\Delta_M(t)= \Delta_M(t^{-1})$ 
and normalized so that $\Delta_M(1)=\left|\Tors\ H_1(M)\right|$.
The RT torsion determines the Alexander polynomial:
$$
\Delta_M(t)= -\hbox{pr}([\tau](M,\xi)) (t-1)(t^{-1}-1) + |\Tors\ H_1(M)|
$$
where pr$:\Z[H_1(M)] \to \Z[t^{\pm 1}]$ 
is induced by the choice of a generator $t$ of $H_1(M)$ up to torsion,
and where $\xi$ is any Euler structure with zero Chern class (see \cite[\S II.5.2]{Turaev_bigbook}). 
So, one gets
$$
\nabla_M(z)  = 1 + |\Tors\ H_1(M)|^{-1} \cdot z^2 
\cdot \hbox{pr}([\tau](M,\xi))|_{z^2=t+t^{-1}-2}.
$$
Observe that the above variable substitution is such that 
$$
z^2=t+t^{-1}-2=-(t-1)(t^{-1}-1) \in I^2.
$$
So, by Corollary \ref{cor:rank_one} applied to $d=2k+3$, the coefficient of $z^{2k}$
in the Alexander--Conway polynomial $\nabla_M(z)$ 
is a rational-valued finite-type invariant of degree $2k$.
This also follows from works by Garoufalidis--Habegger \cite{GH} and Lieberum \cite{Lieberum}.

\subsubsection{The linking pairing}

Recall from \cite[\S X.2]{Turaev_bigbook} that the RT torsion 
determines the linking pairing $\lambda_M: H_1(M) \times H_1(M) \to \Q/\Z$ 
of a rational homology oriented $3$-sphere $M$:
$$
\forall x,y \in H_1(M), \ \tau(M,\xi) \cdot (x-1)\cdot (y-1) = -\lambda_M(x,y)\cdot \Sigma 
\ \  \in \ \frac{\Q[H_1(M)]}{\Z[H_1(M)]}
$$
where $\Sigma:= \sum_{x\in H_1(M)} x$. 
Lemma \ref{lem:general} with $d=1$ implies the basic fact that the  isomorphism in homology
induced by a Torelli surgery preserves the linking pairing.
 
\subsubsection{The cohomology ring}

Turaev has also shown that the leading term of the torsion
with respect to the $I$-adic filtration  can be computed from the cohomology ring.
For instance, take $G=\Z^b$ with $b\geq 3$, and consider the invariant
$$
u: \mathcal{M}(G) \longrightarrow \Hom(\Lambda^3(G^*),\Z), \
(M,\psi) \longmapsto u(M,\psi)
$$
which assigns to $(M,\psi)$ the push-out
by $\psi^*: H^1(M) \to G^*= \Hom(G,\Z)$ of the triple-cup products form
$$
u_M: H^1(M) \times H^1(M) \times H^1(M) \longrightarrow \Z.
$$
It is shown in \cite[\S III.2]{Turaev_bigbook} that, for any Euler structure $\xi$, 
$\tau(M,\xi) \in I(\Z[H_1(M)])^{b-3}$ and 
\begin{equation}
\label{eq:cohomology}
\tau(M,\xi) \ \mod \ I(\Z[H_1(M)])^{b-2} =
\left\{\begin{array}{ll}
q_{b-3}(\hbox{Det}(u_M)) & \hbox{if } b \hbox{ is odd},\\
0 & \hbox{if } b \hbox{ is even}. 
\end{array} \right.
\end{equation}
Here, for $L$ a free Abelian group of rank $b$, 
$\hbox{Det}: \Hom(\Lambda^3(L),\Z) \to \hbox{S}^{b-3}(L^*) $
is a certain polynomial concomitant of degree $b-1$, 
and $q_{k}: \hbox{S}^{k}(L^*) \to I(\Z[L^*])^{k} / I(\Z[L^*])^{k+1}$ 
is the natural  homomorphism defined by 
$g_1 \cdots g_k \mapsto [(g_1-1)\cdots (g_k-1)]$.

Since $u(M,\psi)$ is a degree $1$ finite-type invariant 
(see, for instance, \cite[Lemma 3.2]{Massuyeau_bis}), 
we deduce that $\tau(M,\xi,\psi)\in I^{b-3}/I^{b-2} $ is a finite-type invariant of degree $\leq b-1$
(which does not depend on $\xi$). This also follows from Theorem \ref{th:nonQHS} applied to $d=b$. 

\begin{remark}
Turaev's result (\ref{eq:cohomology}) also reveals
that the degree in Theorem \ref{th:nonQHS} should be improved in many cases (when $b=b_1(M)$ is even).
\end{remark}

\subsection{Domination by finite-type invariants}

\label{subsec:domination}

In view of the results from \S \ref{sec:proofs}, the question of the domination
of the RT torsion by finite-type invariants involves
the nilpotent residue of the augmentation ideal of a group ring.
For finitely generated nilpotent groups, this residue has been computed 
in \cite[Theorem VII.3.1]{Passi}. 
In our situation, their result specializes to the following

\begin{theorem}[Parmenter--Passi--Sehgal \cite{Passi}]
\label{th:residue}
Let $G$ be a finitely generated Abelian group,
whose  torsion subgroup decomposes as
$\oplus_{p\in \mathcal{P}}\ G_p$ into a direct sum of $p$-groups 
indexed by the set of prime numbers $\mathcal{P}$. Then, we have that
$$
\bigcap_{n \geq 0} I(\Z[G])^n = \sum_{\substack{p, q \in \mathcal{P}\\ p\neq q}}
 I(\Z[G_p]) \cdot I(\Z[G_q]) \cdot \Z[G].
$$
\end{theorem}

This theorem together with Lemma \ref{lem:equivalence_tau}
and Lemma \ref{lem:general} imply the following

\begin{corollary} 
\label{cor:domination}
Let $G$ be a finitely generated Abelian group
whose torsion part is reduced to a $p$-group  ($p$ being a prime number).
Then, finite-type invariants of closed oriented connected $3$-manifolds
with Euler structure and homology parametrized by $G$,
dominate the RT torsion.
\end{corollary}

\begin{example}
In  their paper \cite{BL}, Bar-Natan and Lawrence have computed the LMO--Aarhus invariant \cite{LMO,BGRT}
on lens spaces. In particular, they have found that it does not detect the difference between
$M_1:=L(25,4)$ and $M_2:=L(25,9)$.
In contrast, the Reidemeister torsion is well-known to classify lens spaces. 
In particular,  we have that
$$
\tau(M_1,\psi_1,\xi^{\parallel}_1) \neq \tau(M_2,\psi_2,\xi_2^{\parallel}) \in Q(\Z[G]) 
$$
for any homological parametrizations $\psi_1$ and $\psi_2$ by $G:=\Z_{25}$,
and where $\xi_i^{\parallel}$ denotes the unique Euler structure of $M_i$ coming from a parallelization
(i.e., with trivial Chern class). By Corollary \ref{cor:domination},
$(M_1,\psi_1)$ and $(M_2,\psi_2)\in \mathcal{M}(G)$ are distinguished by finite-type invariants for any $\psi_1$ and $\psi_2$, 
although $M_1$ and $M_2$ are not distinguished by rational-valued finite-type invariants 
(as follows from the universality of LMO--Aarhus).
\end{example}

\subsection{The MT torsion}

\label{subsec:MT}

Let $G$ be a finitely generated Abelian group. 
Recall from Example \ref{ex:Milnor--Turaev_torsion} that the MT torsion is this reduction
of the RT torsion  that ignores the torsion part of $G$. It refines the Alexander function or, 
which is essentially the same in dimension three, the Alexander polynomial \cite{Turaev_first}.

Let $\pr: G \to L:= G/ \Tors\ G$ be the canonical projection 
onto the torsion-free quotient of $G$. The \emph{exponential} map of $L$
$$
\exp: L \longrightarrow \widehat{\hbox{S}}(L)\otimes \Q = 
\prod_{ k\geq 0} \hbox{S}^k(L) \otimes \Q ,
\quad l \longmapsto  \sum_{k\geq 0} l^k \otimes 1/k!
$$
takes its values in the completed symmetric algebra of $L$ tensored with $\Q$.
It extends linearly to a ring homomorphism $\exp:\Z[L] \to \widehat{\hbox{S}}(L)\otimes \Q$.

For $(M,\xi,\psi) \in \mathcal{ME}(G)$, we define
$$
T(M,\xi,\psi):= \left\{\begin{array}{ll}
\exp \circ \pr \left(\tau(M,\xi,\psi)\right) & \hbox{ if } \rk(G) > 1,\\
\exp \circ \pr \left([\tau](M,\xi,\psi)\right) & \hbox{ if } \rk(G) = 1,\\
0 & \hbox{ if } \rk(G) = 0.\\
\end{array}\right.
$$
Let also $T_k(M,\xi,\psi)\in \hbox{S}^k(L)\otimes \Q $ be the degree $k$ part of $T(M,\xi,\psi)$.

\begin{corollary}
\label{cor:Milnor--Turaev_torsion}
For closed oriented connected $3$-manifolds with Euler structure and homology parametrized by $G$, 
the invariant $T=\prod_{k\geq 0} T_k$  is equivalent to the MT torsion. 
Moreover, for each $k\geq 0$, $T_k$ is a finite-type invariant of degree $\leq k+2$. 
In particular, the MT torsion is dominated by finite-type invariants.
\end{corollary}

\begin{proof}
One easily checks that the map 
$$
\varphi_k: L\longrightarrow \hbox{S}^k(L), \ l \longmapsto l^k
$$ 
is \emph{polynomial} of degree $\leq k$, in the sense that its linear extension
to $\Z[L]$ vanishes on $I(\Z[L])^{k+1}$. See \cite[Chapter V]{Passi}. 
It follows from Theorem \ref{th:nonQHS} and Corollary \ref{cor:rank_one} that
$$
T_k(M,\xi,\psi)= \left\{\begin{array}{ll}
\varphi_k \circ \pr\left(\tau(M,\xi,\psi)\right) \otimes 1/k! & \hbox{ if } \rk(G) > 1\\
\varphi_k \circ \pr \left([\tau](M,\xi,\psi)\right) \otimes 1/k! & \hbox{ if } \rk(G) = 1\\
\end{array}\right.
$$
is a finite-type invariant of degree $\leq k+2$. 

Moreover, the group homomorphism $\varphi_k: I(\Z[L])^k/I(\Z[L])^{k+1} \to \hbox{S}^k(L)$
induced by $\varphi_k$ can easily be computed:
$$
\forall l_1, \dots, l_k \in L, \ \
 \varphi_k(\left[(l_1-1) \cdots (l_k-1)\right])= k! \cdot l_1 \cdots l_k.
$$
One the other hand, there is this homomorphism 
$q_k: \hbox{S}^k(L) \to I(\Z[L])^k/I(\Z[L])^{k+1}$ 
defined by $l_1 \cdots l_k \mapsto [ (l_1-1) \cdots (l_k-1) ]$. It is surjective and, 
since $\varphi_k \circ q_k = k! \cdot \Id$, it is injective as well. So, $\varphi_k$ is injective too.
Since the ideal $I(\Z[L])$ is residually nilpotent, the map $\exp: \Z[L] \to \widehat{\hbox{S}}(L)\otimes \Q$ is injective. 
We conclude that the invariant $T$ is equivalent to the MT torsion.
\end{proof}

\begin{remark}
These finiteness properties of the MT torsion parallel
those of the  Alexander polynomial of oriented links in $\mathbf{S}^3$,
which have been shown by H. Murakami in \cite{Murakami}.
\end{remark}

\begin{remark} 
If one restricts $\xi$ to be an Euler structure $\xi^\parallel$ coming from a parallelization, the quantity
$T(M,\xi,\psi)$ does not depend  anymore on $\xi$. So, the invariant
$$
T: \mathcal{M}(G) \longrightarrow \widehat{\hbox{S}}(G/ \Tors\ G)\otimes \Q,  \quad
(M,\psi) \longmapsto T(M,\xi^\parallel,\psi)
$$
is available. In particular, its degree $k$ part 
$$
T_k(M,\psi) \in  \hbox{S}^k(L)\otimes \Q \simeq 
\hbox{S}^k(L\otimes \Q) 
$$
defines a symmetric $k$-multilinear form on 
$$
\Hom_\Q(L\otimes \Q,\Q) \simeq  
\xymatrix{{\Hom_\Z(G,\Q)} & \ar[l]_-{\psi^*}^-\simeq {H^1(M;\Q).}}
$$
That form is easily seen to coincide (up to multiplication by $(-2)^k$) 
with Turaev's \emph{$m$-th moments} form of the torsion function \cite[\S IX.2]{Turaev_bigbook}.
Thus, the latter is a finite-type invariant of degree $\leq k+2$.
\end{remark}

\bibliographystyle{amsalpha}

\vspace{1cm}

\noindent
\textsc{\footnotesize Dipartimento di Matematica Applicata,
Via Bonanno Pisano 25/B, 56126 Pisa, Italia}\\
\emph{\footnotesize E-mail:}
\texttt{\footnotesize massuyeau@mail.dm.unipi.it}\\

\vspace{0.5cm}

\noindent
Current address:\\[0.2cm]
\noindent
\textsc{\footnotesize Institut de Recherche Math\'ematique Avanc\'ee
(CNRS/Universit\'e Louis Pasteur),  7 rue Ren\'e Descartes,
67084 Strasbourg, France}\\
\emph{\footnotesize E-mail:}
\texttt{\footnotesize  massuyeau@math.u-strasbg.fr}

\end{document}